\documentclass[reqno,12pt]{amsart}
\usepackage{amssymb,graphicx}
 \def\Z{{\mathbb Z} } \def\R{{\mathbb R} }
\def\jj{{\bf j} } \def\ii{{\bf i} }
\def\ss{{\bf s} } \def\tt{{\bf t} }
\def\aa{{\rm a} } \def\bb{{\rm b} } \def\cc{{\rm c} }
\def\AA{$\mathcal A$\ } \def\BB{$\mathcal B$\ } \def\CC{$\mathcal C$\ }
\def\DD{$\mathcal D$\ } \def\EE{$\mathcal E$\ } \def\FF{$\mathcal F$\ }
\def\GG{$\mathcal G$\ }
\def\inte{{\rm int}\, }

\newtheorem{Thm}{Theorem} [section]
\newtheorem{Def}[Thm]{Definition}
\newtheorem{Prop}[Thm]{Proposition}

\newtheorem{Rem}[Thm]{Remark}
\newtheorem{Exa}[Thm]{Example}
\newtheorem{Prob}[Thm]{Problem}
\newtheorem{Conj}[Thm]{Conjecture}
\newenvironment{Proof}[0]
{\medskip \noindent {\bf Proof.} \ }{\ \hfill $\Box$\medskip}

\begin{document}
\title{Combinatorial topology of three-dimensional
self-affine tiles}
\author{Christoph Bandt}
%\address{Institute of Mathematics, University of Greifswald, Germany}
\date{\today}
\begin{abstract}
We develop tools to study the topology and geometry of self-affine
fractals in dimension three and higher. We use the self-affine
structure and obtain rather detailed information about the
connectedness of interior and boundary sets, and on the dimensions
and intersections of boundary sets. As an application, we describe
in algebraic terms the polyhedral structure of the six fractal
three-dimensional twindragons. Only two of them can be
homeomorphic to a ball but even these have faces which are not
homeomorphic to a disk.
\end{abstract}
\maketitle

\section{Introduction}
When L\'evy \cite{le} introduced his famous curve in 1938, he also
constructed fractal surfaces in a similar way. 70 years later, we
have plenty of papers on fractal sets in the plane, and a number
of general statements and constructions which hold in every
dimension, but very few studies on the geometry of fractals in
dimension $\ge 3.$

\begin{figure}
 \begin{center}
  \includegraphics[width = 150mm]{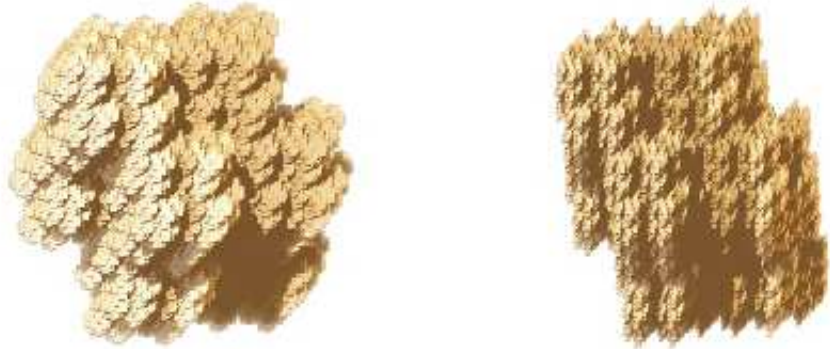}
  \caption{\label{twin6} Two views of the twindragon \CC , made by chaoscope \cite{de}.
  Is this set homeomorphic to a ball? Is the interior connected?
  We develop concepts and tools to answer such questions.}
  \end{center}
\end{figure}

Visualization is certainly one reason to prefer plane sets: they
can be easily shown on a computer screen. In dimension 3,
visualization is more difficult, even though there are good
ray-tracing programs, like chaoscope \cite{de} and IFS builder
\cite{km}, which allow to look at three-dimensional fractals from
any chosen viewpoint, with prescribed light sources. Figure
\ref{twin6} shows two views of one of our main examples. Does it
have interior points? If so, is the interior connected? Is the
figure simply connected? Unlike in two dimensions, the answer to
such questions can hardly be guessed from the pictures - this
requires tools which we develop here.

There are several mathematical difficulties in dimension $\ge 3.$
We avoid the problem that linear maps do usually not commute, by
concentrating ourselves to fractals generated by a single matrix.
But there is another problem which will concern us: proper
self-similarity is rather an exceptional case. The moduli of all
eigenvalues of a $3\times 3$ matrix do rarely coincide while for a
$2\times 2$ matrix they coincide as soon as the eigenvalues are
complex. This problem is addressed in Section 2, and Theorem 2.3
implies that in ${\mathbb R}^3$ there is essentially only one type
of `integer' self-similar tile with 2 up to 7 pieces.

Other difficulties concern the topology of three-dimensional
fractals. Our main aim here is to study the topology of
self-affine tiles, and of their boundaries which are fractal
surfaces. In particular, we would like to know which of these
tiles are homeomorphic to a closed ball. In two dimensions, the
corresponding question for disk-like tiles was discussed in a
large number of papers \cite{at1,at2,bg,bw,gg,lrt,ll,st,t} (more
references can be found in \cite{ll,v,w}), and Jordan curve
arguments were heavily used. Among others, it is easy to see that
connectedness of the interior of a self-affine tile in the plane
is sufficient for its disk-likeness. In ${\mathbb R}^3$ this is
not the case. The topology of the boundary of a tile can be much
more complicated than in ${\mathbb R}^2.$ Even the connectedness
of a tile is not so easy to verify \cite{ag,klr}.

The geometric phenomena which we encounter in fractal tiles are
not caused by knots or wild spheres \cite{bc}. It is rather the
self-affine fibre structure of the boundary which makes the
topology complicated. Apparently, the different eigenvalues of the
generating matrix produce long and thin fibres which can pierce
the interior of neighboring tiles, or distort the boundary
structure.

We develop algebraic tools which describe the geometry of a
self-affine tile $T$ in arbitrary dimension, in a similar way as
homotopy and homology groups describe the geometry of manifolds.
Our tools use self-affinity in a specific way, and they yield
quite detailed information about the geometry. The basic concept
is the {\it neighbor graph} $G=(V,E).$ It can be considered as a
blueprint of $T$ which contains all information about the topology
of $T$ in a simple scheme. In the case where $T$ tiles by
translations in a unique way, the vertices $v\in V$ are those
translations for which $T+v$ belongs to the tiling and intersects
$T.$ Such translations $v$ will be called neighbor maps, and they
do not only represent the neighboring piece $T+v$ but also the
boundary set $T\cap (T+v)$ of $T.$ This set can be a face, an
`edge' or just a point of $T,$ or a more complicated fractal
boundary set. For a ball-like tile we expect that faces are
homeomorphic to disks and edges are homeomorphic to intervals, and
we have to check this idea.

Neighbors in lattice tilings were introduced by Indlekofer,
K\'atai and Racsko \cite{ikr}, and neighbor maps for self-similar
sets were defined by Bandt and Graf \cite{bg}. For a self-affine
tile, the boundary sets have themselves a self-affine structure:
each face or edge is the union of smaller copies of other boundary
sets. This is indicated by the edges of the graph $G.$ The
definition of $G$ in Section 3 implies that the boundary $\partial
T$ is a self-affine graph-directed construction in the sense of
Mauldin and Williams \cite{mw}, and the graph which directs this
construction is $G$ itself. This fact was stated by Scheicher and
Thuswaldner \cite{st} for lattice tiles, and used in many other
papers to determine the Hausdorff dimension and topological
structure of $\partial T$ for plane self-similar tiles
\cite{gi,sw,klsw,lrt,ll,t,ve}. Some authors
\cite{dkv,st,at1,at2,al} prefer the contact matrix defined by
Gr\"ochenig and Haas \cite{gh} which describes a subgraph of $G.$

A first application of the neighbor graph concerns the existence
of self-affine lattice tilings, a problem treated in \cite{gh} and
in a series of papers by Vince \cite{v1,v4,v}. Theorem 5.1 shows
that the existence of such tilings is equivalent to the
compatibility of all neighbors, which can be easily checked with
$G.$

In sections 7 to 12 we study the topology of the boundary of $T.$
In section 8 we define the hierarchy of different boundary sets
like faces, edges, points etc. which corresponds to the hierarchy
of sets with different Hausdorff dimensions in Mauldin-Williams
constructions. An important problem is to determine all {\it
faces} -- those boundary sets which cover an open set on the
boundary. Two methods are developed for this purpose. Theorem 9.2
gives an algorithmic approach which only uses the combinatorics of
the neighbor graph. Theorem 10.1 provides an analytic approach
based on recent work by He and Lau \cite{hl} and Akiyama and
Loridant \cite{al}, connected with eigenvalues of adjacency
matrices and Hausdorff dimension. In section 11 we define a
polyhedral structure of fractal tiles by the intersection sets of
two and more faces. There are corresponding neighbor intersection
graphs $G^2,G^3,...$ which allow to calculate details of the
polyhedral structure and compare it with the structure of ordinary
polyhedra.

In contrast to the contact matrix, all these tools can also be
applied to tiles with incompatible neighbors as well as to
fractals which are not tiles, and neighbors can be related by
arbitrary isometries of ${\mathbb R}^3$ instead of translations.
To illustrate our methods, we selected a small class of examples,
which is introduced in section 6: the {\it twindragons.} These are
the self-affine lattice tiles with two pieces. The one-dimensional
twindragon is the interval. In the plane, there are three
examples: rectangle or parallelogram, the twindragon and the tame
twindragon. They are all disk-like \cite{ikr,bg}. In
three-dimensional space, there are seven examples, which we denote
by letters \AA  to \GG . While \AA is conjugate to a cube, \FF and
\GG  have an extremely intricate structure. We mention some of
their properties without proof. For \DD  and \EE  we determine the
boundary faces and their intersections, and we show that their
interior is not simply connected. \BB  and \CC  seem to be
homeomorphic to a ball. In section 12 we determine their exact
polyhedral structure. We show that they have faces which are not
homeomorphic to a disk, but \CC has connected interior.

This paper focusses on concepts rather than algorithmic questions
although it seems clear that the neighbor graph and related finite
automata are well-suited for computerized evaluation. Section
\ref{algo} gives a short account of algorithmic aspects.

\section{Self-affine lattice tiles}
{\bf Basic concepts. } The following definitions are standard
\cite{lw,lw1,v,w}. A {\it lattice} $L$ in $\R^n$ is the set of
integer linear combinations of $n$ linearly independent vectors
$e_1,...,e_n.$ Usually we take $L=\Z^n.$ A regular-closed set $T$
is called a {\it tile} if it admits a tiling of $\R^n.$ A {\it
tiling} by the tile $T$ is a countable union $\bigcup_k T_k$ which
covers $\R^n,$ such that each $T_k$ is isometric to $T,$ and no
two tiles have interior points in common. We have a {\it lattice
tiling} by $L$ if the tiles $T_k$ are just the translates $T+k$
with $k\in L.$

A tile $T$ is a {\it self-affine lattice tile} if there is an
affine expanding mapping $g$ and a lattice $L$ such that $g$
preserves $L$ and maps $T$ to a union of tiles $T+k_i$ with
$k_i\in L.$ With respect to the basis vectors $e_i,$ the map $g$
has a matrix representation $g(x)=Mx$ where $M$ is an integer
matrix. Expanding means that all eigenvalues of $M$ have modulus
greater one. $T$ is determined by the map $g$ and by the lattice
coordinates $k_1,...,k_m$ of the tiles which form the supertile
$g(T).$
\begin{equation} g(T)=MT=\bigcup_{j=1}^m T+k_j \label{selfa}
\end{equation}
For $T$ to become a tile, we must have $m=|\det M|,$ and the set
$K=\{k_1,...,k_m\}$ must fulfil some condition. A sufficient
condition is \cite{b5}
\begin{equation} g(L)=\bigcup_{j=1}^m L+k_j=L+K \label{selfb}
\end{equation}
in which case $K$ is called a {\it standard digit set} or {\it
complete residue system.} The latter name comes from the fact that
for $|\det g|=m$, the subgroup $g(L)$ of the additive group $L$
has $m$ residue classes, so by (\ref{selfb}) each class is
represented by exactly one $k_i.$ See Lagarias and Wang
\cite{lw,lw1} who also studied non-standard digit sets. In this
paper, we consider only standard digit sets.

The self-affine tile $T$ is called {\it self-similar} if the
mapping $g$ is a similarity mapping with respect to the Euclidean
metric. In this case, we must distinguish the standard basis which
defines the Euclidean distance, and the basis for the lattice $L$
which need not be orthonormal.

{\bf Conjugacy of tiles. } $T$ is said to be {\it conjugate to a
self-similar tile} if there is a linear map $h$ so that
$\tilde{T}=h(T)$ is a self-similar tile. The expanding map for
$\tilde{T}$ is $\tilde{g}=hgh^{-1},$ the lattice is
$\tilde{L}=h(L),$ and the digits are $\tilde{k}_i=h(k_i),$ for
which (\ref{selfa}) and (\ref{selfb}) are easily checked.

\begin{Exa} In the plane, we consider the digits
$k_1={0\choose 0}, k_2={1\choose 0}.$ The parallelogram $T$ with
vertices ${0\choose 0}, {1\choose 0}, {1\, \choose 1}, {2\choose
1}$ is a self-similar tile with respect to $M={1\ 1\choose 1\,
-1}$ and the lattice $L=\Z^2.$ The rectangle $T'=[0,1]\times
[0,\sqrt{2}]$ is also a self-similar lattice tile, with $M'={0\
\sqrt{2} \choose \sqrt{2}\ 0}$ and $L'=\{ {m\choose\sqrt{2}\, n}\,
|\, m,n\in\Z\} .$ The unit square $T''=[0,1]^2$ is a self-affine
tile with respect to the matrix $M''={0\, 2 \choose 1\, 0}$ and
$L''=\Z^2,$ but not self-similar. However, $T''$ is conjugate to
the self-similar tile $T'$ since $h(x)={1\,\  0 \choose
0\,\sqrt{2}}x$ fulfils $h(T'')=T'.$ Note that $\tilde{h}(T)=T'$
for $\tilde{h}(x)={1\, -1 \choose 0\,\sqrt{2}}x.$ \end{Exa}

We shall identify conjugate tiles so that $T,T',T''$ are
considered as different versions of the same tile, which is
essentially self-similar. {\it Considering tiles up to conjugacy
means that we focus on the essential data: the characteristic
polynomial of $g$ instead of different matrix representations.} We
shall also assume that $k_1={0\choose 0}$ since this can be
obtained by conjugacy with a translation.

\begin{Prop} (The class of self-similar tiles)\\
A self-affine lattice tile $T$ with respect to $g(x)=Mx$ is
conjugate to a self-similar tile if and only if all eigenvalues of
$M$ have the same modulus.
\end{Prop}

\begin{Proof} Necessity of the condition follows from the fact
that eigenvalues are not changed under conjugacy. On the other
hand, if the eigenvalues have equal modulus, there is a matrix $B$
such that $\tilde{g}(x)=BMB^{-1}x$ is a similarity map. Thus
$h(x)=Bx$ maps $T$ to a self-similar tile.\end{Proof}

\begin{Thm} (Very few self-similar lattice tiles in dimension 3)\\
If in $\R^3$ a self-affine lattice tile $T$ with respect
to $g(x)=Mx$ is conjugate to a self-similar tile, then one of the
following two conditions is true.
\begin{enumerate}\item[ (i)] $m=|\det
M|$ is a cubic number - in particular $m\ge 8,$
\item[(ii)] $M$ is conjugate to $\tilde{M}=\left(\begin{array}{lcr}0&0&\pm m\\
1&0&0\\ 0&1&0\end{array}\right) .$
\end{enumerate}
\end{Thm}

\begin{Proof} With respect to the lattice base, $g$ has an integer
matrix, and so the characteristic polynomial $p(\lambda)$ of $M$
has integer coefficients. We express them in terms of the
eigenvalues $\lambda_i.$
\[ p(\lambda)=-\lambda^3 +(\lambda_1+\lambda_2+\lambda_3)\cdot\lambda^2
-(\lambda_1\lambda_2+\lambda_1\lambda_3+\lambda_2\lambda_3)\cdot\lambda
+\lambda_1\lambda_2\lambda_3 \] where
$\lambda_1\lambda_2\lambda_3=\det M=\pm m.$ Since $T$ is conjugate
to a self-similar tile, all $\lambda_i$ have equal modulus:
$|\lambda_i|=r>1.$

If the eigenvalues are real, then $\lambda_i=\pm r$ and $m=\pm
r^3.$ In this case $\lambda_1+\lambda_2+\lambda_3\in \{\pm r, \pm
2r, \pm 3r\},$  so $r$ must be an integer and $m$ a cubic number.

Now let us assume $\lambda_1,\lambda_2$ are complex eigenvalues.
Then $\lambda_1=r(\cos\alpha +i\sin\alpha)$ and
$\lambda_2=r(\cos\alpha -i\sin\alpha)$ for some $\alpha ,$ so
\[ p(\lambda)=-\lambda^3 +2r(\cos\alpha +\frac{s}{2})\cdot\lambda^2
-(r^2+2sr^2\cos\alpha)\cdot\lambda +r^3s
\] where $s=\pm 1$ is the sign of $\lambda_3.$ The coefficient of
$\lambda$ can be written as $-2sr^2(\cos\alpha +\frac{s}{2}).$
Since the coefficients are integers, either $r$ must be a rational
number - and hence an integer, and $m$ a cubic number. Or
$\cos\alpha +\frac{s}{2}$ must be zero. In this case,
$p(\lambda)=-\lambda^3+sm$ which is the characteristic function of
$\tilde{M}.$
\end{Proof}

\begin{Rem}
The orthogonal part of $\tilde{M}$ is a rotation around $120^o$
for $s=+1,$ and a rotation around $60^o$ composed with a
reflection at the plane of rotation for $s=-1.$ Even in the case
where $r$ is an integer and $m$ a cubic number, $2r\cos\alpha$
must be an integer and $\alpha$ can assume only few values. Thus
there are really very few self-similar lattice tiles in $\R^3.$
\end{Rem}

\section{Neighbors of tiles and self-affine sets}
When we want to study the boundary $B$ of a lattice tile $T,$ it
is quite natural to consider the {\it neighbors} in the tiling -
those tiles $T+k$ for which $B_k=T\cap (T+k)$ is non-empty. The
$B_k$ can be considered as the faces of $T,$ and $B$ is the union
of the $B_k.$ Obviously, the number of neighbors is finite. This
idea was introduced by Gilbert \cite{gi} and Indlekofer, K\'atai
and Racsko \cite{ikr} and used in many other papers.
\smallskip

{\bf Neighbors in fractals. } There is a related concept of
neighbor in self-similar fractals which appeared first in
\cite{bg}. Suppose we have a fractal $A$ which consists of two
copies or itself: $A=f_1(A)\cup f_2(A),$ where the $f_i$ denote
contracting similarity maps. Then the geometry and topology is
determined by the structure of the intersection set $C=f_1(A)\cap
f_2(A).$  Only at points of $C$ it makes sense to zoom into the
picture. If $x$ is in $f_1(A),$ say, and $U$ is a neighborhood of
$x$ in $A$ which does not intersect $C,$ then the magnified copy
$f_1^{-1}(U)$ of $U$ does already exist in $A.$

The sets $f_1^{-1}(C)$ and $f_2^{-1}(C)$ can be considered as
boundary sets of $A,$ where $A$ intersects a {\it potential
neighbor} $f_1^{-1}f_2(A)$ or $f_2^{-1}f_1(A).$ These need not be
the only boundary sets, however, since $C$ consists of
intersections of smaller pieces like $f_1f_2(A)\cap f_2f_2(A),$
and these may be at other positions when they are pulled back to
$A.$ This argument shows that the boundary sets have a
self-similar structure. The number of boundary sets obtained by
going to smaller and smaller pieces need not be finite, but in
most common examples it is.

Let us define neighbors in this general sense. Let $f_1,...,f_m$
denote contractive affine mappings on $\R^n.$ There is a unique
compact non-empty set $A$ with
\begin{equation}A=\bigcup_{j=1}^m
f_j(A),\label{selfc}\end{equation} the {\it self-affine set} with
respect to the $f_i,$ cf. \cite{fa}. For each integer $q$ the set
$A$ splits into $m^q$ small copies $f_{\jj}(A)=A_{\jj},$ where
$\jj=j_1...j_q$ denotes a word of length $q$ from the alphabet
$I=\{ 1,...,m\},$ and $f_\jj=f_{j_1}\cdot ...\cdot f_{j_q}.$

A potential {\it neighbor of $A$} has the form
$h(A)=f_{\ii}^{-1}f_{\jj}(A)$ where $\ii,\jj\in I^*$ are any words
over $I$ such that the pieces $A_{\ii}=f_{\ii}(A)$ and
$A_{\jj}=f_{\jj}(A)$ intersect. The idea is that $f_{\ii}^{-1}$
maps $A_{\ii}$ to $A$ and ${A_\jj}$ to $h(A),$ a neighbor set
which intersects $A$ and should be of comparable size. $h$ is
called a {\it neighbor map} and the set $B=A\cap h(A)$ the
corresponding {\it boundary set} of $A.$ To get `comparable size',
we confine ourselves to words $\ii,\jj$ of equal length and make
some assumptions concerning the $f_i.$

\begin{Prop} (Neighbor maps which are isometries)\\
Let $A\subset \R^n$ be a self-affine set
of the form (\ref{selfc}). Let $\ii,\jj$ be words of the same
length $q,$ and $h=f_{\ii}^{-1}f_{\jj}.$
\begin{enumerate} \item[ (i)] If all $f_j$ are similarity maps with the
same similarity factor $r,$ then $h$ is a Euclidean isometry.
\item[(ii)] If  $f_j(x)=M^{-1}(x +k_j)$ where $M$ is an expanding matrix
and $k_1,...,k_m\in\R^n,$ then $h$ is a translation.
\end{enumerate}
\end{Prop}

\begin{Proof} (i): If the $f_j$ are similarity maps with factor
$r$ then $f_{\jj}$ and $f_{\ii}^{-1}$ are similarity maps with
factor $r^q$ and $r^{-q},$ respectively. Thus $h$ is an
isometry.\\
(ii) is first proved for $q=1.$ Since $f_i^{-1}(y)=My-k_i,$ we
have
\begin{equation}f_i^{-1}f_j(x)=x+k_j-k_i \label{neibeg}\end{equation}
Moreover, if $h(x)=x+v$ is a translation then
\begin{equation}
f_i^{-1}hf_j(x)=x+Mv+k_j-k_i \label{induc}
\end{equation}
Now induction on $q$ shows that for $h=f_{\ii}^{-1}f_{\jj},$
\begin{equation} h(x)=f_{i_q}^{-1}...f_{i_1}^{-1}
f_{j_1}...f_{j_q}(x)=x+M^{q-1}(k_{j_1}-k_{i_1})+M^{q-2}(k_{j_2}-k_{i_2})
+...+(k_{j_q}-k_{i_q}) \label{neivec}
\end{equation}
\end{Proof}

{\bf Lattice tiles as a special case. } The case (ii) includes our
self-affine lattice tiles since the defining equation
(\ref{selfa}) can be rewritten as
\begin{equation} T=\bigcup_{j=1}^m g^{-1}(T)+g^{-1}(k_j)
=\bigcup_{j=1}^m M^{-1}(T+k_j) \label{selfd}
\end{equation}
We get $f_j(x)=M^{-1}(x+k_j)$ as contracting maps for the
self-affine set $T.$

At this point it is possible to explain the concept of standard
digit set $k_1,...,k_m$ which says that $k_{j_q}-k_{i_q}$ is not
in $ML$ for $j_q\not= i_q$ (see section 2). This property implies
that the translation vector in (\ref{neivec}) can never be zero if
$j_q\not= i_q$ since $k_{j_q}-k_{i_q}$ cannot cancel with the
other terms in (\ref{neivec}). As a consequence, the inequality
$h\not= id$ then holds whenever $\ii\not=\jj .$ From the results
in \cite{bg} then follows the so-called open set condition which
says that the pieces $f_{\ii}(T)$ and $f_{\jj}(T)$ have no
interior points in common. (In \cite{b5}, this was shown by
Baire's category theorem.) In other words, $T$ and a translate
$h(T)$ in (\ref{neivec}) have no common interior points.

\section{The neighbor graph}
{\bf Self-similarity of boundary sets. } To explain the
self-similar structure of boundary sets, we return to the general
case. Consider intersecting pieces $A_{\ii}$ and $A_{\jj}$ in the
self-affine set $A.$ Since all pieces divide into $m$ subpieces,
like $A$ in (\ref{selfc}),
\[ A_{\ii}\cap A_{\jj} =\bigcup_{i,j=1}^m A_{\ii i}\cap A_{\jj j}
.\] Now consider the boundary sets $B=f_{\ii}^{-1}(A_{\ii}\cap
A_{\jj})$ and $B_{ij}=f_{\ii i}^{-1}(A_{\ii i}\cap A_{\jj j}).$
Then
\[ B=f_{\ii}^{-1}\left(\bigcup_{i,j=1}^m A_{\ii i}\cap A_{\jj j}\right)
=\bigcup_{i,j=1}^m f_if_i^{-1}f_{\ii}^{-1}(A_{\ii i}\cap A_{\jj
j}) =\bigcup_{i,j=1}^m f_i (B_{ij}) . \] At $B,B_{ij}$ we
suppressed the subscripts $\ii,\jj$ since such an equation holds
for each possible boundary set. In other words, {\it the
subdivision of pieces induces self-similar representations of the
type} (\ref{selfc}) {\it for all boundary sets} - not with one
type of set, but with several types. We must assume now that we
have only finitely many possible boundary sets, which is true for
lattice tiles with standard digit sets. The unions on the right
contain a lot of empty terms, so we introduce a graph which better
describes the system of equations.\smallskip

{\bf Concept and properties of neighbor graph. } The {\it neighbor
graph} $G=(V,E)$ of a self-affine fractal $A$ has as vertex set
all neighbor maps $h=f_{\ii}^{-1}f_{\jj},$ and a directed edge
marked with $i$ goes from $h$ to
$h'=f_i^{-1}f_{\ii}^{-1}f_{\jj}f_j.$ Loops and multiple edges
(even with the same label) are possible \cite{bm}. For a formal
definition, we focus on self-affine tiles.

\begin{Def}\label{defng}
Let $T=\bigcup f_j(T)$ be a self-affine lattice tile with
$f_j(x)=M^{-1}(x+k_j)$ so that the neighbor maps have the form
$h(x)=x+k.$ We identify $h$ with the translation vector $k,$ and
asopciate it with the boundary set $B_k=T\cap (T+k).$ An edge from
$k$ to $k'$ with label $i$ is drawn when relation (\ref{induc})
holds:
\[ G=(V,E) \quad\mbox{ with }\quad
V=\{ k\, |\, T\cap (T+k)\not=\emptyset\} \qquad \mbox{ and }\]
\begin{equation} E=\bigcup_{i=1}^m E_i\quad\mbox{ with }\quad
E_i=\{ (k,k',i)\, |\, k'=Mk+k_j-k_i \}
\label{neigra}
\end{equation}\end{Def}

Here $E_i$ denotes the set of edges with label $i.$ It should be
mentioned that $V$ contains a root vertex $k=0$ (or $h_0=id$ in
the more general notation) which does not correspond to a boundary
set. There are edges $(0, k_j-k_i,i)$ from the root which
correspond to the first step (\ref{neibeg}) in the calculation of
neighbor maps. The loops $(0,0,i)$ will not be drawn, however. For
the calculation of $G,$ see section \ref{algo}, \cite{bm,br} and
the example below.

\begin{Prop} (The boundary equations)\\
Let $T=\bigcup_{j=1}^m f_j(T)$ be a self-affine lattice tile with
neighbor graph $G=(V,E).$ The boundary sets $B_k$ corresponding to
the $k\in V\setminus \{0\}$ fulfil the following equations.
\begin{equation}
B_k=\bigcup\, \left\{ f_j(B_{k'}\, |\, j\in\{ 1,...,m\} ,\,
(k,k',j)\in E\right\} \label{gradir}\end{equation}
\end{Prop}

This is an immediate consequence of Definition \ref{defng} and of
the discussion above. Mauldin and Williams \cite{mw} called
families of such fractals $B_k$ {\it graph-directed
constructions}, and proved that the $B_k$ are uniquely determined
by the graph $G$ and the maps $f_j.$ At the end of section 10 we
briefly discuss the open set condition of (\ref{gradir}).

\begin{Prop} (The basic symmetry of $G$)\\
If $h$ is a neighbor map, then $h^{-1}$ is also a neighbor map.
For translation this means that to every $k$ there is a $-k,$ and
each boundary set $B_k=T\cap (T+k)$ has an opposite set
\[ B_{-k}=(T-k)\cap T=B_k -k.\] This symmetry of $V$ will extend as
a graph isomorphism to the edges, but not to the labels. If
$k'=Mk+k_j-k_i$ then $(k,k')$ has label $i$ and $(-k,-k')$ has
label $j.$
\end{Prop}

In a previous publication with M. Mesing \cite{bm}, we wrote both
labels $i,j$ at the edge $(k,k').$ This is not necessary here
since we can always recover the second label from the opposite
edge.

{\bf A two-dimensional example. } In our examples, neighbors will
be denoted by lower-case Roman letters, and $-\bb$ will denote the
opposite vertex of b. We start with a two-dimensional tile
\cite{b5,v} which can be considered as a modification of the
square, and as an extension of the Sierpi\'nski gasket. The
construction of neighbor graphs for three-dimensional tiles
proceeds in the same way as demonstrated here.

\begin{figure}
 \begin{center}
  \includegraphics[width = 150mm]{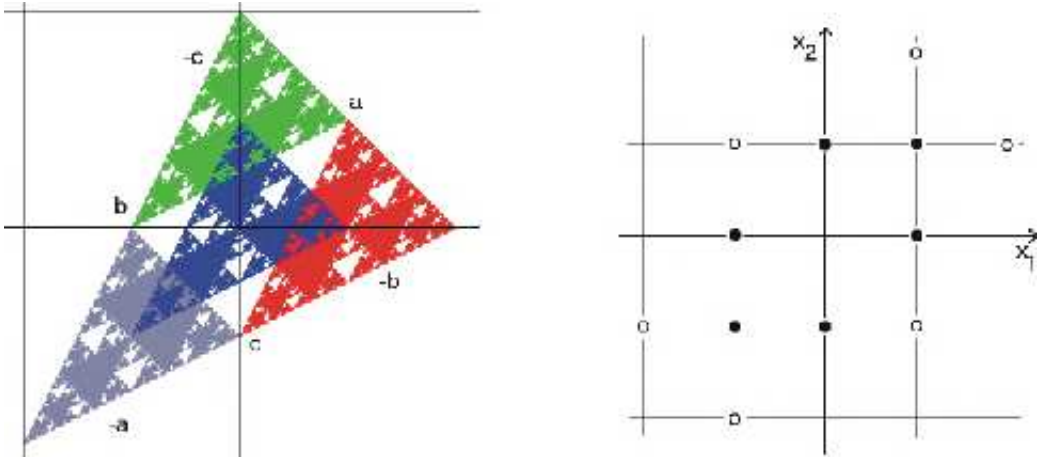}
  \caption{\label{siertile} A two-dimensional tile
  and its neighbor translations.}
  \end{center}
\end{figure}

\begin{Exa}
Let $f_j(x)=(x+k_j)/2$ for $j=1,...,4$ where $k_1={0\choose 0},
k_2={1\choose 0}, k_3={0\choose 1},$ and $k_4={-1\choose -1}.$ The
tile $T$ is shown in Figure \ref{siertile}. Clearly,  $T\subset
[-1,+1]^2.$ Thus a translation $x={x_1\choose x_2}$ with $T\cap
(T+x)\not=\emptyset$ must fulfil $|x_1|\le 2$ and $|x_2|\le 2.$
Moreover, $L$ is the integer lattice, so $x_1,x_2\in\{
-2,-1,0,1,2\} .$ In this way we prove that for an arbitrary
self-affine lattice tile, the set $V$ is finite.

The tile $T$ has six neighbors which intersect $T$ in a single
point, and six neighbors which intersect $T$ in an uncountable set
which is in fact a Sierpi\'nski gasket. \end{Exa}

\begin{Proof} We start with the translation vectors $k=k_j-k_i$
where $i\not= j,$ see the right-hand part of Figure
\ref{siertile}. If $k$ is in $V,$ there is the edge $(0,k,i).$
Since our matrix is $M={2\ 0\choose 0\ 2},$ all these points do
belong to $V,$ with a loop $(k,k,j)$ at each point $k,$ because
$k_j-k_i=2(k_j-k_i)+k_i-k_j.$

Moreover, these will be the only points of $V.$ For the other
points $x$ with $|x_1|\le 2, |x_2|\le 2$ the recursion
(\ref{induc}) will lead to vectors $k'$ which are outside the
range of $x_1,x_2$ and thus are no neighbor maps. A simple
calculation \cite{bm} shows that once we are outside the range, we
can never come in again by the recursion (\ref{induc}).

Thus the translations $\pm\aa=\pm {1\choose 1}, \pm\bb=\pm
{-1\choose 0}, \pm\cc=\pm {0\choose -1}$ marked by $\bullet ,$ and
$\pm {2\choose 1}, \pm {1\choose 2}, \pm {1\choose -1}$ marked by
$\circ$ are the vertices of our neighbor graph. The last six
vertices have only the loop and no further outgoing edges. These
six vectors are the translations along the sides of the big
triangle in Figure \ref{siertile}. They translate $T$ to a
neighbor $h(T)$ which has a single point as intersection with $T.$
These neighbors are called {\it point neighbors.}

To find point neighbors, we need only look at the graph, not at
the picture. The loop at $k={2\choose 1}$ with label $2$ says that
the boundary set has address $2222...=\overline{2},$ and the
opposite vector ${-2\choose -1}$ has address $\overline{4}.$ This
means that $T\cap (T+k)$ is exactly one point, which has address
$\overline{2}$ in $T$ and address $\overline{4}$ in $T+k.$ Indeed
the intersection is the point $k_2={1\choose 0}$ which is the
fixed point of $f_2,$ and $T+k$ will meet this point with the
translate of $k_4,$ the fixed point of $f_4.$ The translate of
$k_3$ along the vector $k'={-1\choose 1}$ is also $k_2.$ This way
it turns out that two of the six point neighbors will meet $T$ at
each vertex of the triangle.

\begin{figure}
 \begin{center}
  \includegraphics[width = 60mm]{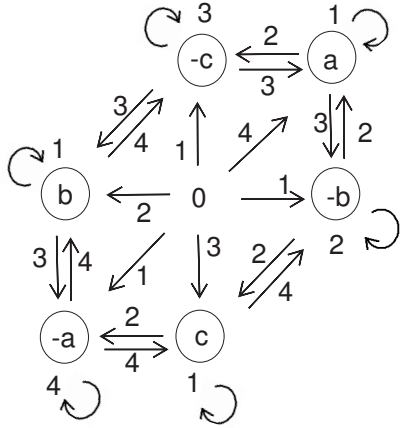}
  \caption{\label{siergra} The reduced neighbor graph for
  the Sierpi\'nski tile.}
  \end{center}
\end{figure}

We now forget the point neighbors $\circ .$  The graph of the
remaining six $\bullet$ neighbors $\pm\aa , \pm\bb, \pm\cc$ and
the root, shown in Figure \ref{siergra}, will be called the
reduced neighbor graph \cite{bm}. On the left of Figure
\ref{siertile}, the corresponding boundary sets are indicated. As
mentioned above, each vertex $k=k_j-k_i$ has a loop with label
$i.$ The edge $(-\bb,\aa)$ has label 2 and the opposite edge has
label 3 because $\aa=2\cdot(-\bb)+k_3-k_2.$ In a similar way, all
other edges were determined.

The graph in Figure \ref{siertile} is irreducible which means that
all corresponding boundary sets have a similar structure. This
also follows from the equations of the graph-directed system
(\ref{gradir}): \[ B_\aa =f_1(B_\aa )\cup f_2(B_{-\cc})\cup
f_3(B_{-\bb}), \ B_{-\bb}=f_2(B_\aa\cup B_{-\bb}\cup B_\cc), \
B_{-\cc}=f_3(B_\aa\cup B_\bb\cup B_{-\cc}) .
\] All these boundary sets are in fact small Sierpi\'nski gaskets,
as can be seen in the figure. The theorem below shows that they
are uncountable.
\end{Proof}

{\bf Topology and address map. } A sequence $\ss =s_1s_2...\in
I^\infty$ is called address of the point $x$ in a self-affine tile
$T$ if it is contained in the pieces $T_{s_1...s_n}$ for
$n=1,2,...$ (cf. \cite{bar}, Chapter IV, \cite{fa}, Section 9.1).
The points of the intersection sets $T_i\cap T_j$ have two
addresses, one starting with $i$ and the other one with $j.$ The
address map $\pi: I^\infty\to T$ is a quotient map, it determines
the topology of $T.$ For that reason it is important to know which
addresses will be identified. The neighbor graph provides this
information, as can be easily shown (\cite{bm}, sections 4 and 5).
A path in the directed graph $G$ is a finite or infinite sequence
$e_1e_2...$ such that the terminal vertex of $e_k$ coincides with
the initial vertex of $e_{k+1}.$ A finite path $e_1...e_n$ is a
cycle if the terminal vertex of $e_n$ is the initial vertex of
$e_1.$

\begin{Thm}(Neighbor graph and addresses of points, see \cite{bm})\\
Let $T$ be a self-affine tile with the address map $\pi:
I^\infty\to T$ and the neighbor graph $G=(V,E).$
\begin{enumerate} \item[ (i)] Two sequences $\ss=s_1s_2...$ and
$\tt=t_1t_2...$ are addresses of the same point if and only if
there is a path $e_1e_2...$ starting in the root with labels
$s_1s_2...,$ such that the opposite path -$e_1$-$e_2...$ has
labels $t_1t_2...$
\item[(ii)] The addresses of the points of a boundary set $B_k$
coincide with label sequences of paths starting in the vertex $k.$
They form a regular language $L_k.$
\item[(iii)] If only one cycle can be reached from the vertex $k$
by a directed path, then $B_k$ is a singleton.
\item[(iv)] If two different cycles can be reached from $k,$ and
these cycles can be reached from each other, then the set $B_k$ is
uncountable. \end{enumerate}
\end{Thm}
Note that a cycle with a diagonal, or with a loop at one of its
points counts as two different cycles which can reach each other.

\section{Existence of self-affine lattice tilings}
Before we can go on, we must settle one difficulty with neighbors.
They need not be compatible: $h(T)$ and $\tilde{h}(T)$ can be
neighbors of $T,$ but not of each other. In that case, different
tiles in a tiling will have different neighborhoods.

It turns out that the compatibility of neighbors is connected with
another question which was studied by Gr\"ochenig and Haas
\cite{gh} and in several papers of Vince \cite{v1,v4,v} as a
``central'' problem in the field. A self-affine lattice tile with
standard digit set will always admit self-affine tilings. It is a
non-trivial result that it will also admit lattice tilings
\cite{lw2} but these need not be self-affine in the sense that
tiles assemble to form supertiles (\cite{v}, Example 4.2). Vince
(\cite{v4},\cite{v}, Theorem 4.3) gave 10 equivalent conditions
for the existence of a self-affine lattice tiling. Two of them are
of an algorithmic nature.

Here we extend Vince's list by proving that a self-affine lattice
tiling exists if and only if all neighbors are compatible with
each other, and we show that this property can be easily decided
with the neighbor graph. The following theorem also shows that all
self-affine sets with mappings $f_j(x)=M^{-1}(x+k_j)$ which are
not lattice tiles must have incompatible neighbors.

\begin{Thm} (Compatible neighbors and self-affine lattice tilings)\\
Let $A\subset \R^n$ be a connected self-affine set with respect to
mappings $f_j(x)=M^{-1}(x+k_j), j=1,...,m$ where $M$ is an
expanding matrix and the open set condition holds. Then the
following conditions are equivalent.
\begin{enumerate}
\item[ (i)] All potential neighbors appear together at one piece $A_{\ii}.$
\item[(ii)] $A$ is a self-affine lattice tile which admits a
self-affine tiling by a lattice.
\item[(iii)] In the neighbor graph, each vertex $k\not= 0$ has at least
one incoming edge with each of the labels $j=1,...,m.$
\end{enumerate}\end{Thm}

\begin{Proof} (i)$\Rightarrow$(ii): We assume $A$ is a self-affine set, and
$A_{\ii}$ has all possible neighbors. Then the number of neighbors
must be finite, at most $m^n$ if $\ii=i_1...i_n.$ Moreover,
$f_{\ii}^{-1}(A)$ consists of copies of $A$ which include all
possible neighbors $h(A).$ Since the maps are translations
$h(x)=x+k$ by Proposition 3.1, let again $K$ denote the set of all
these translation vectors. We define an infinite pattern of
translates of $A$ as a union of increasing compact sets: \[
\bigcup_{q=1}^\infty f_{\ii}^{-q}(A) .\] All translates $A'$ of
$A$ in this pattern have the form $f_{\ii}^{-q}(A_{\ii\jj})$ for
some number $q$ and some word $\jj .$ This implies that they all
have the same neighbors as $A.$ To see this, consider the
(potential) boundary of a subpiece $A_{\ii j}$ with $j\in\{
1,...,m\},$ defined as the union of all its intersections with
potential neighbors. This boundary is contained in the union of
the boundary of $A_{\ii}$ and the intersections $A_{\ii j}\cap
A_{\ii i}$ with $i=1,...,m, i\not= j.$ So by assumption the
boundary of $A_{\ii j}$ is contained in $A,$ and by induction on
the length of $\jj$ this is proved for subpieces $A_{\ii\jj}.$

Thus all `atoms' $A'$ of our infinite pattern have a maximal set
of neighbors which covers their boundary completely. But there is
only one maximal set: the complete set of translates $A'+k$ with
$k\in K.$

Let $L$ denote the lattice generated by $K.$ Since $A$ was
connected, any two pieces of the same level are connected by a
chain of neighboring pieces (cf. \cite{br}, 8.2.1), and each atom
has the form $A'=A+k$ with $k\in L.$ Moreover, $L$ has the same
rank as the linear subspace generated by $A,$ which we will now
assume is $\R^n.$
% This is true because our infinite self-affine pattern must extend
% infinitely  into all directions of the linear space, and it is
% contained in $A+L$ where $A$ is compact.

Since all $A'$ have the neighbor vectors $K,$ the infinite pattern
coincides with $A+L.$ As a consequence, the pattern must be a
tiling of $\R^n:$ If there was a small open set $U$ outside $A+L$
then there would be arbitrary large copies $f_{\ii}^{-q}(U)$
outside $A+L$ which would contradict the finite mesh size of $L.$
Thus we have constructed a self-affine lattice tiling by copies of
$A.$
\smallskip

(ii)$\Rightarrow$(iii): Let $\R^n=A+L$ be a self-affine lattice
tiling by the tile $A.$ Then each tile of the tiling has the same
set of neighbor translations $K.$  Thus if $k$ is in $K,$ then $k$
appears as a neighbor translation of any piece $T_j$ of a
supertile $T,$ where $j\in\{ 1,...,m\}.$ If the neighbor is
another piece $T_i$ of $T,$ there is the edge $(0,k,j)$ in $G.$ If
the neighbor of $T_j$ is in another supertile $T'$ then we have
the edge $(k',k,j)$ where $k'$ is the neighbor map between $T$ and
$T'.$ Thus for each vertex $k$ of the vertex set $K,$ there is an
incoming edge with label $j.$\smallskip

(iii)$\Rightarrow$(i): We assume that each vertex $k$ in the
neighbor graph $G=(V,E)$ has incoming edges with each label $j.$
We have to find a word $\ii=i_1...i_N$ such that $A_\ii$ has all
neighbors $k\in V.$ For each $k$ there must be a suffix
$i_{n(k)}...i_N$ of $\ii$ which consists of labels of a path from
$0$ to $k.$

It is not difficult to construct such a word $\ii$ by induction,
going paths {\it backwards} (against the direction of edges)
towards 0. Start with a vertex $k_1$ and let $\jj_1=j_1...j_{q_1}$
be the labels of a backward path from $k_1$ to 0. According to our
assumption, there are also backward paths with label sequence
$\jj_1$ from all other vertices. Their endpoints, which are
different from 0, will form a set $V_1.$

Next, take a $k_2\in V_1$ and a backward path
$\jj_2=j_{q_1+1}...j_{q_2}$ from $k_2$ to 0. Also take backward
paths from all other $k\in V_1$ with label sequence $\jj_2$ and
denote the set of their non-zero endpoints by $V_2.$ Continue with
$k_3\in V_2,$ and so on.  Since $V_{n+1}$ has less points than
$V_n,$ we will have $V_n=\emptyset$ for some $n.$ The sequence
$\ii=j_{q_n}...j_{q_{n-1}}...j_{q_1}...j_1$ will have all required
suffixes.
\end{Proof}

All examples considered in this paper will fulfil the condition of
theorem 5.1. For Example 4.4 this can be seen in Figure
\ref{siergra}, for the twindragons in the figures below.

\section{The seven three-dimensional twindragons}
{\bf Twindragons and their symmetry. } We now define the family of
examples we are going to study here. A twindragon is a self-affine
lattice tile with $m=2$ pieces. The expanding matrix $M$ has
determinant $\pm 2.$ In dimension 1, a twindragon is an interval.
In $\R^2,$ there are three examples: the rectangle (Example 2.1),
the ordinary twindragon and the tame twindragon, and they are all
disk-like \cite{bg}.

Since we identified conjugate tiles, we could take $k_1=0$ in
section 2. For twindragons in $\R^n,$ we can also choose $k_2$
arbitrarily, as long as $k_2$ is not in an invariant linear
subspace of $M.$ For the basis $\{ k_2, g(k_2),...,g^{n-1}(k_2)\}$
the matrix $M$ of $g$ will always be the same, determined by the
characteristic polynomial of $g$ (see \cite{gg1} and the proof of
Theorem 6.2 below). For our twindragons in $\R^3,$ we choose
another normalization:
\begin{equation}\label{knorm} k_1=(-\frac12 ,
0,0)', \quad k_2=(\frac12 ,0,0)' \,
\end{equation}
This has the effect that 0 is the symmetry center of $T,$ and in
the neighbor graph, an edge with label 1 passes from the root to
$\aa =(1,0,0)'.$

\begin{Prop} (Symmetry of twindragons)\\
Each twindragon $T$ in $\R^n$ has a symmetry center at
$c=\frac12(k_1+k_2).$ The point reflection at $c$ interchanges the
two pieces $T_1,T_2,$ and each boundary set $B_k$ with the
opposite set $B_{-k}.$ Moreover, each boundary set $B_k$ of $T$
also has a symmetry center at $c_k=c+\frac{k}{2}.$
\end{Prop}

\begin{Proof} We can take $k_1,k_2$ from (\ref{knorm}) since
the symmetry is not changed by an affine conjugacy. Thus
$f_1(x)=M^{-1}x+k_1, f_2(x)=M^{-1}x-k_1.$ With $\phi(x)=-x$ we
obtain $\phi f_1=f_2\phi$ or $f_i=\phi f_{3-i}\phi^{-1}$ since
$\phi^{-1}=\phi.$ For any word $\ii=i_1...i_n$ on $I=\{ 1,2\},$ we
have $f_\ii=\phi f_{\jj}\phi^{-1},$ where each $j_q=3-i_q$ is the
opposite symbol. Since the point of $T$ with address
$\ss=i_1i_2...$ can be represented as $x=\lim_{n\to\infty}
f_{i_1}f_{i_2}...f_{i_n}(z)$ where $z$ is any starting point
\cite{fa}, this implies that $\phi$ transforms this point into the
point with the opposite address. In particular $\phi(T)=T$ and
$\phi(T_i)=T_{3-i}$ from which it follows that $\phi$ maps
$T_1\cap T_2$ onto itself.

For each boundary set $B_k=T\cap (T+k)$, the map $\phi$ transforms
$T+k$ to $\phi(T)-k=T-k,$ thus $\phi (B_k)=B_{-k}.$ On the other
hand we had $B_{-k}=B_k-k,$ so $\psi(x)=-x+k$ maps $B_k$ to
itself. $\psi$ is the point reflection at $c_k=\frac{k}{2}.$
\end{Proof}

{\bf The seven twindragons. } The last proposition shows that the
tile generated by $f_1,f_2$ can also be generated by $f_1\phi$ and
$f_2\phi,$ where $\phi$ is the point reflection at $c.$ Since in
$\R^3$ the map $\phi$ has determinant -1, this has the consequence
that we need only consider the matrices $M$ with $\det M=2.$ With
-2 we get the same twindragons.

\pagebreak

\begin{figure}
 \begin{center}
  \includegraphics[width = 150mm]{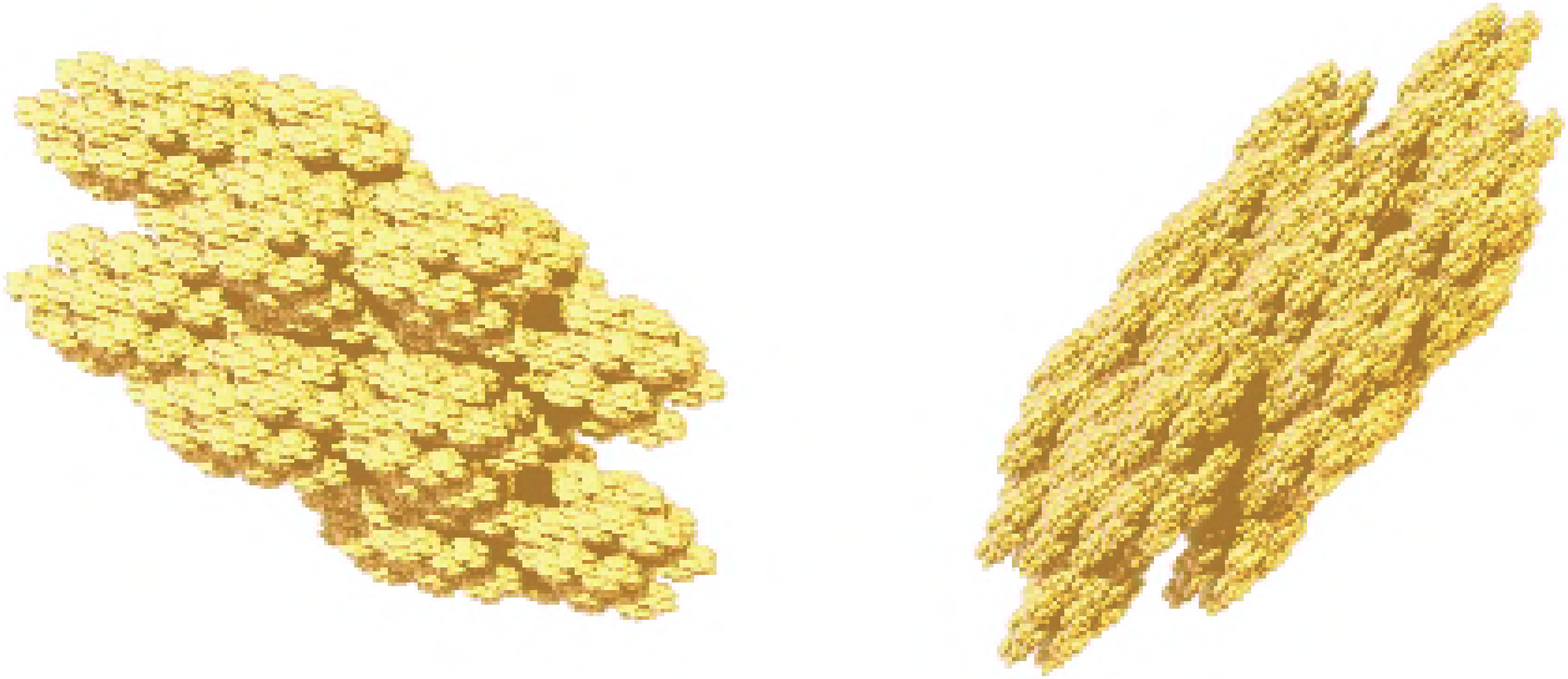}
  \caption{\label{twin7} Twindragon \BB, shown first in \cite{gg1}}
  \end{center}
\end{figure}

\begin{Thm} (List of three-dimensional twindragons, cf.
\cite{klr,ag})\\
Up to conjugacy, there are exactly seven twindragons, each of
which is uniquely determined by the pair of coefficients $(a,b)$
of its characteristic polynomial
\[ p(\lambda)=-\lambda^3+a\lambda^2+b\lambda+2\, :\]
\begin{center}\begin{tabular}{|l|c|c|c|c|c|c|c|}\hline
Twindragon &\AA&\BB&\CC&\DD&\EE&\FF&\GG\\ \hline
Parameter&$(0,0)$&(-1,1)&(1,-1)&$(0,1)$&(2,-2)&$(1,0)$&$ (0,2)$\\
\hline \end{tabular}
\end{center}
\end{Thm}

\begin{Proof}
These are the only integers $(a,b)$ for which $p(\lambda)$ has
only roots of modulus $>1$ \cite{klr,ag}. For $k_1=0,
k_2=(1,0,0)'$ and the basis $\{ k_2, Mk_2, M^2k_2\}$ the affine
map
\[ g(x)=\left(\begin{array}{lcr}0&0&2\\
1&0&b\\ 0&1&a\end{array}\right)\cdot x \] has characteristic
polynomial $p(\lambda ),$ and yields a twindragon. The
standardization (\ref{knorm}) translates this twindragon by
$(-\frac12 , 0,0)'$ so that its center is 0.
\end{Proof}

\begin{Exa} ( The non-fractal twindragon \AA )\\
\AA , the self-affine cube with 2 pieces, is the three-dimensional
analogue of Example 2.1. When considered as a rectangular
parallelepiped with side lengths $1, \sqrt[3]{2},\sqrt[3]{4},$ it
is the self-similar exception mentioned in Theorem 2.3, (ii). The
similarity map $g$ then is a $120^o$ rotation composed with a
homothety by the factor $\sqrt[3]{2}.$ Of course, this tile has 26
neighbors: 6 faces, 12 edge neighbors and 8 point neighbors.
\end{Exa}

\begin{figure}
 \begin{center}
  \includegraphics[width = 150mm]{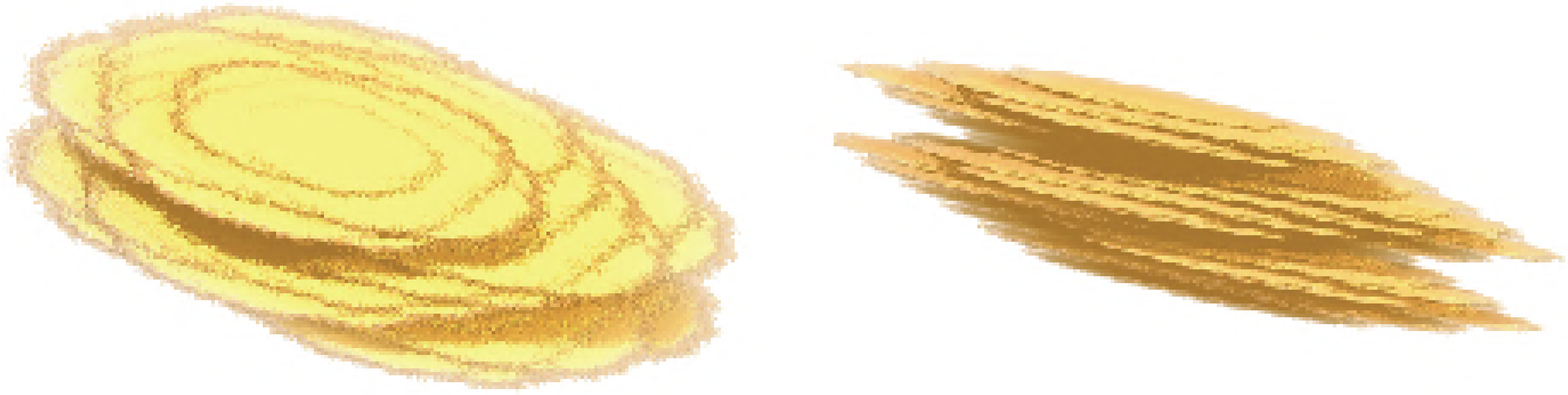}
  \caption{\label{twin2} Twindragon \EE }
  \end{center}
\end{figure}

{\bf Number of neighbors of twindragons. } The twindragons form a
nice little family with seven examples on which we can test our
neighbor methods. However, the six fractal examples are truly
self-affine, which makes their topology intricate. As a first
overview, we compare the moduli of eigenvalues of $M$ - the
inverses of the contraction factors of the $f_j$ along different
directions - and the number of neighbors. If the contraction
factors differ much, there are lots of neighbors. We also
determine the number of infinite boundary sets, by deleting from
the neighbor graph the vertices representing finite boundary sets,
as was explained in Example 4.4 and Theorem 4.5. It turns out that
all infinite boundary sets are uncountable. The number of faces
will be determined in the sections below.
\smallskip

\begin{Prop} (First neighbor calculation for the seven twindragons)
\begin{center}\begin{tabular}{|l|c|c|c|c|c|c|c|}\hline Twindragon &\AA&\BB&\CC&\DD&\EE&\FF&\GG\\ \hline
Parameter&(0,0)&(-1,1)&(1,-1)&(0,1)&(2,-2)&(1,0)&(0,2)\\
\hline  Complex eigenvalues'
modulus&1.26&1.29&1.22&1.15&1.14&1.09&1.06\\ \hline
Real eigenvalue&1.26&1.21&1.35&1.52&1.54&1.70&1.77\\
\hline Number of neighbors&26&18&20&34&34&48&76\\ \hline
Uncountable boundary sets&18&14&12&14&32&48&76\\ \hline
Faces&6&14&12&14&12&16&$\ge 22$\\ \hline
\end{tabular}\end{center}
\end{Prop}\smallskip

For \BB , \CC and \DD the number of infinite boundary sets is
smaller than for the cube, so that we can expect them to be more
or less ball-like. \EE has medium complexity, and the last two
examples have many infinite neighbor sets -- they could be
three-dimensional counterparts of Example 4.4.

When we compare with the figures, we see that the twindragons
possess a kind of fibre structure along the direction of the
smaller eigenvalues of $M,$ that is, the larger eigenvalues of
$M^{-1}.$ For \BB , the pieces drawn in Figure \ref{twin7} are
stretched along the direction of the real eigenvalue of $M^{-1}$
which is the largest contraction factor. So roughly the pieces
look like `cigars'. In all other examples, the larger contraction
factor is given by the complex eigenvalues, and the pieces look
like leaves or plates. For \CC in Figure \ref{twin6}, the leaves
still have some thickness, while for \EE illustrated in Figure
\ref{twin2}, the plates are already very thin. For the other three
twindragons, the plates look similar or still thinner.\smallskip

{\bf How to draw neighbor graphs of twindragons. } Let us mention
some other nice properties of our family and fix some notation for
the neighbor graphs. First, all twindragons are connected. For two
pieces this follows from $T_1\cap T_2\not=\emptyset$ (see Barnsley
\cite{bar}, 8.2.1). If this relation would not hold, $T$ would be
a Cantor set and would have no interior points.

Next, the only possible differences $k_j-k_i$ are
$k_2-k_1=(1,0,0)', k_1-k_2=(-1,0,0)',$ and $k_i-k_i=0.$ The
corresponding labels $i$ of the edges are 1 and 2 and $i$ where
the last case means two edges labelled with 1 and 2. We shall draw
edges for the last case as double arrows, and we distinguish
labels 1 and 2 by assigning a fat tip to the arrows with label 2.
This way we save the labelling of edges. The label of the opposite
edge is $j=3-i$ for simple arrows, and $j=i$ for double arrows.

In view of the symmetry, we can simplify the neighbor graph
further by drawing only one vertex from each pair of opposite
vertices $k,-k$ and introducing the convention that a wavy arrow
$\leadsto$ to $k$ denotes an arrow to $-k.$ This will also
simplify the equation systems (\ref{gradir}) for the boundary
sets. It really matters whether one has to work with 20 or 40
vertices! Note that the opposite vertex of $k$ is $-k$ not only as
a vector, but by Proposition 6.1 also as a boundary set:
$B_{-k}=-B_k.$

For sake of brevity, we shall provide no details for the
complicated twindragons \FF and \GG , but we give full arguments
for the other four examples.

\begin{figure}
 \begin{center}
\hfill\begin{minipage}[b]{60mm}\includegraphics[width =
55mm]{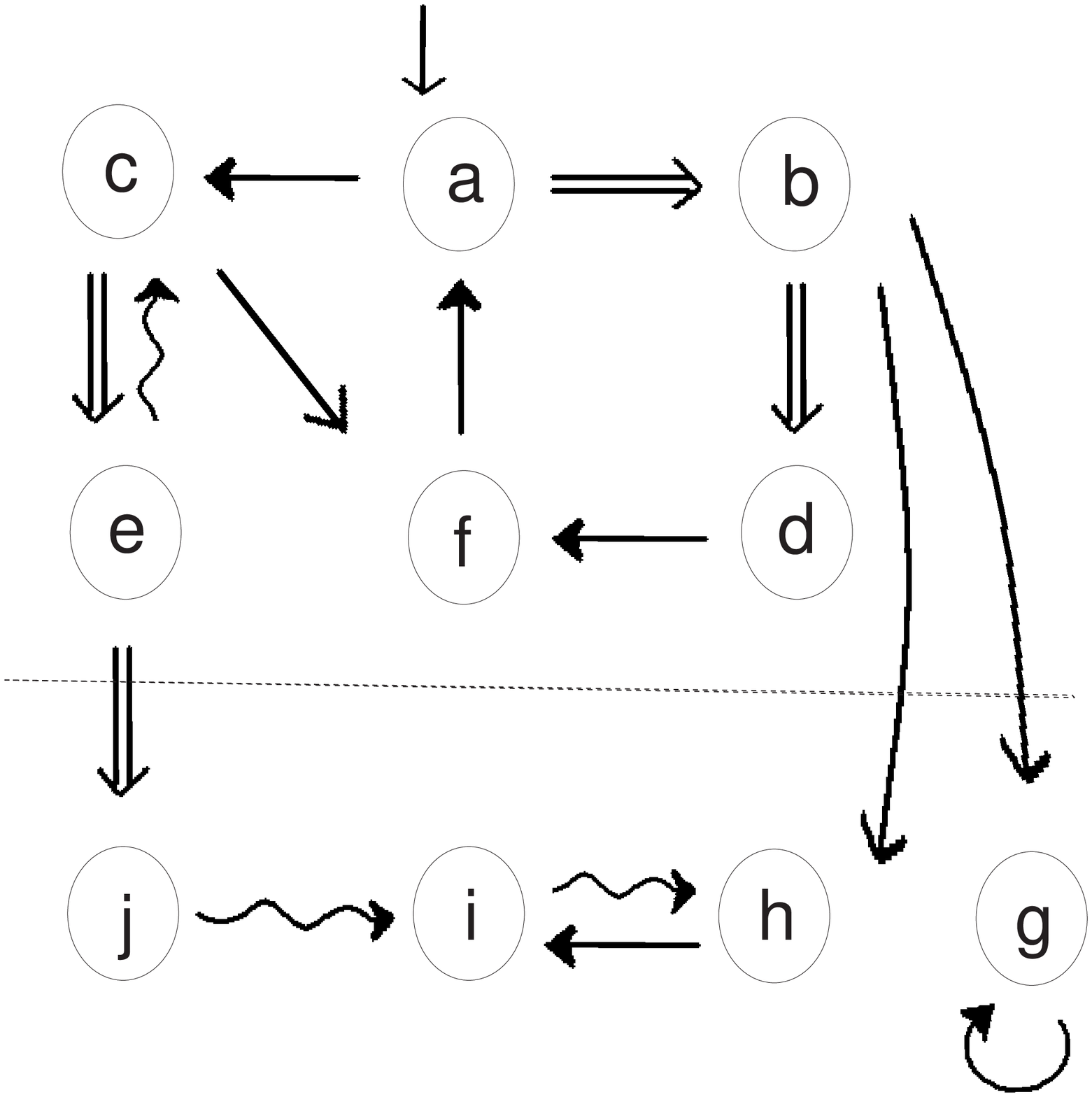}
\end{minipage}
\hfill
\begin{minipage}[b]{70mm}\begin{tabular}{|c|c|c|c|c|c|c|c|c|c|}\hline
a&b&c&d&e&f&g&h&i&j\\ \hline 1&0&-1&0&0&1&1&-1&1&2\\ \hline
0&1&1&0&-1&-1&0&0&-2&-1\\ \hline 0&0&0&1&1&1&1&1&1&0\\ \hline
\end{tabular}\end{minipage}
  \caption{\label{ng6} The neighbor graph for \CC contains all information
  about the topology of Figure \ref{twin6}. The translation
  vectors, which we shall not need,
  indicate the position of the center of the neighbor.}
  \end{center}
\end{figure}

\begin{figure}
 \begin{center}
  \includegraphics[width = 65mm]{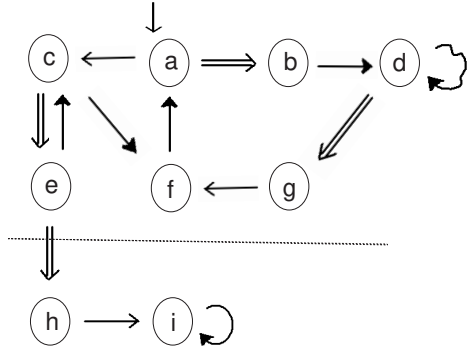}
  \caption{\label{ng7} Neighbor graph for
  the twindragon \BB. There are two point neighbors: i with address
  $\overline{2}$ and h with address $1\overline{2}$.}
  \end{center}
\end{figure}

To verify the condition of Theorem 5.1, we check in Figures
\ref{ng6}, \ref{ng7}, \ref{ng4}, and \ref{ng2} that each vertex is
reached either by a double arrow, or by two ordinary arrows with
different tips indicating label 1 and 2, or by an ordinary and a
wavy arrow with equal tips. In the last case the wavy arrow leads
to $-k,$ so we have to count the opposite edge.

\begin{Exa} We derive the equations (\ref{gradir}) of the boundary sets from
Figure \ref{ng6}. We neglect all point neighbors, given by
vertices below the thin line.
\[ B_\aa=f_1(B_\bb)\cup f_2(B_\bb)\cup f_2(B_\cc),\quad
B_\bb=f_1(B_{\rm d})\cup f_2(B_{\rm d}), \quad B_{\rm
d}=f_2(B_{\rm f}),\]
\[ B_\cc=f_1(B_{\rm e})\cup f_2(B_{\rm e})\cup
f_1(B_{\rm f}),\quad B_{\rm e}=f_2(-B_\cc ), \quad B_{\rm
f}=f_2(B_\aa )\, .\] By substitution we reduce the system to two
equations:
\[ B_\aa=f_2(B_\cc)\cup\bigcup_{i,j=1}^{2} f_{ij22}(B_\aa)\, ,\qquad
B_\cc=f_{12}(B_\aa)\cup\bigcup_{i=1}^{2} f_{i2}(-B_\cc)\, .\]
\end{Exa}\smallskip

\section{Identifying faces: the simple case}
{\bf Defining faces. } The topological boundary of a set $T$ will
be denoted by $\partial T= T\setminus\inte T.$ If $T$ is a tile in
a tiling, then $\partial T$ must be covered by the boundary of
neighboring tiles: $\partial T=\bigcup_{k\in V} B_k.$

Theorem 4.5 says that we can determine the cardinality of a
boundary set $B_k$ by counting the number of paths starting from
$k$ in the neighbor graph $G.$  It is much more complicated to
determine the topology of $B_k.$ In three-dimensional self-affine
tiles, we would like to distinguish boundary sets which are
two-dimensional fractals, one-dimensional curves, Cantor sets or
just finite sets.

Let us say that the boundary set $B_k=T\cap (T+k)$ is a {\it face}
of $T$ if there is a point $x\in B_k$ and an open neighborhood $U$
of $x$ with $U\cap\partial T\subset B_k$ and $U\subset T\cup
(T+k).$ In that case, all boundary points of $U\cap \inte T$ must
belong to $B_k.$ Thus $B_k$ has topological dimension 2 since the
boundary of each open set in $\R^n$ has dimension $n-1.$ The
following statement shows that we have six faces in Example 4.4.

\begin{Prop} Let $k\not= 0$ be a vertex in $G$ such that for each
other vertex $k'\not= 0,$ either $k'$ or the opposite vertex $-k'$
can be reached by a path from $k.$ Then $B_k$ is a face of $T.$
\end{Prop}

\begin{Proof} Each tile must have faces. By (\ref{gradir}),
any path from $k$ to $k'$ implies the existence of a small copy of
$B_{k'}$ in $B_k.$ Thus our assumption implies that $B_k$ contains
copies of faces, so it must be a face (cf. Proposition 8.1).
\end{Proof}

\begin{figure}
 \begin{center}
  \includegraphics[width = 55mm]{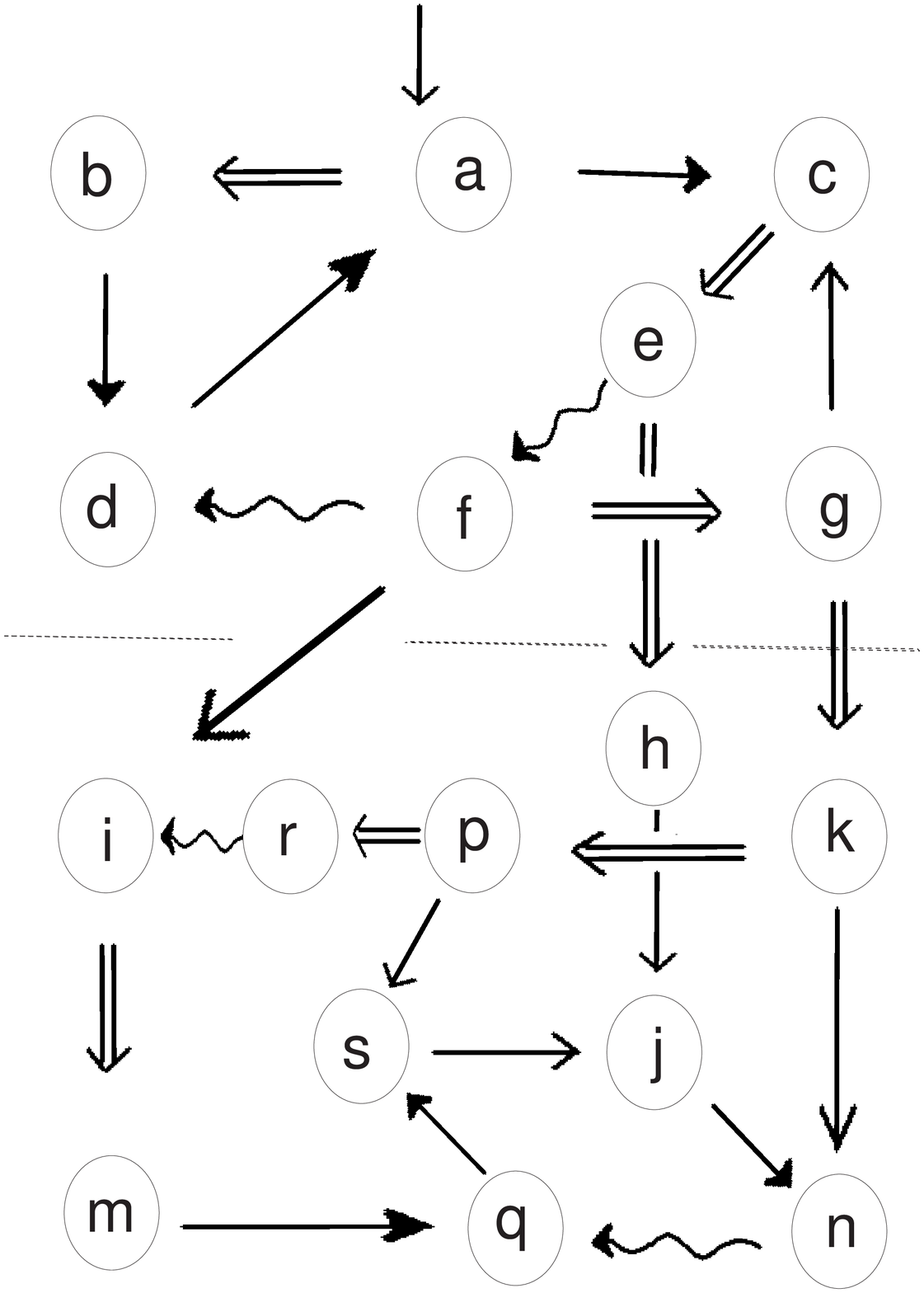}
  \caption{\label{ng4} Neighbor graph for
  the twindragon \DD. The neighbors i,r,p and k describe
  finite boundary sets with more than one point.}
  \end{center}
\end{figure}

For twindragons, there is always an edge with label 1 from the
root to the first other vertex $\aa=k_2-k_1=(1,0,0)'.$ and this
must be a face by Proposition 7.1. It turns out that for the
twindragons \BB, \CC and \DD, all infinite boundary sets fulfil
the assumption of the proposition since from their vertices in
$G,$ there is a path to $\aa .$ The neighbor graphs are shown in
Figures \ref{ng6}, \ref{ng7} and \ref{ng4}. The point neighbors
are separated from the faces by a thin line. They are found by
Theorem 4.5, (iii).

For \BB this can be seen in Figure \ref{ng7}. For \CC we have the
loop at g with address $\overline{2},$ and the 2-cycle $\{{\rm i}
, {\rm h}\}$ in Figure \ref{ng6} which denotes a 4-cycle $\{ {\rm
i, -h, -i, h}\}$ in the complete neighbor graph, with address
$\overline{2112}.$ (The wavy arrow indicates that we go to the
opposite sets, and there the labels 1 and 2 have to be
interchanged. A second wavy arrow reestablishes the original
labels.) In Figure \ref{ng4} we have a similar 4-cycle $\{ {\rm
q,s,j,n}\}$ which denotes an 8-cycle of point neighbors in the
complete neighbor graph with address $\overline{21221211}.$ Also h
and m are point neighbors, while i,r,p and k describe finite
boundary sets with more than one point.  We summarize:

\begin{Prop} (Faces of \BB, \CC and \DD)\\
The twindragon \CC has 12 faces and 8 point neighbors.
\BB has 14 faces and 4 point neighbors. The twindragon \DD has 14
faces, 12 point neighbors and 8 further finite boundary sets.
\end{Prop}

\section{Different types of boundary sets}
{\bf The partial order for boundary sets. } To compare different
infinite boundary sets, we define an order on the vertices of $G.$
For $u,v\in V$ we write $u\succ v$ if there is a directed path
from $u$ to $v.$ We say that the vertices are equivalent, and
write $u\sim v$ if either $u=v$ or both $u\succ v$ and $v\succ u.$
For the boundary sets, $u\succ v$ means that $B_u\supset f_\ii
(B_v)$ for some word $\ii=i_1...i_q.$

\begin{Prop} If $u\succ v$ then $\dim B_u\ge \dim B_v$ holds for
the topological dimension, the Hausdorff dimension and the box
dimension. Moreover, if $B_v$ is a face, then also $B_u$ will be a
face.
\end{Prop}

\begin{Proof} The dimension concepts mentioned here (and also some
others \cite{fa}) have the following properties.
\begin{enumerate}
\item[ (i)] If $A\subseteq B$ then $\dim A \le \dim B .$
\item[(ii)] If $h$ is an affine map then $\dim h(B)=\dim B.$
\end{enumerate}
Together with the remark above this proves the first assertion.
\\ If $B_v=T\cap T+v$ is a face, there is $x\in B_v$ with an open
neighborhood $U\subset T\cup T+v .$ Now $u\succ v$ means that
$B_u\supset f_\ii (B_v)$ for some word $\ii$ where $T+u\supset
f_\ii (T+v).$ So $y=f_\ii (x)$ and its neighborhood $f_\ii (U)$
can be taken to show that $B_u$ is a face.
\end{Proof}

{\bf Strong components of the neighbor graph. } We can now
classify the vertices of $G$ in a similar way as states of a
Markov chain. Clearly $\sim$ is an equivalence relation. The
topological and Hausdorff dimension, as well as the property of
describing a face, are invariant with respect to $\sim .$ Let
$\tilde{V}$ be the set of equivalence classes with respect to
$\sim$ which are called strong components of $G.$ Then we get a
new graph $\tilde{G}$ by drawing an edge from
$(\tilde{u},\tilde{v})$ if there is $u\in\tilde{u}$ and
$v\in\tilde{v}$ such that $(u,v)$ is an edge in $G.$

The graph $\tilde{G}$ compresses the structure of $G.$ If two
vertices $u,v\in G$ belong to the same class, this means that
$B_u$ contains a copy of $B_v$ and conversely, so both boundary
sets have essentially the same structure. The root $0\in V$ would
also be the root of $\tilde{G}$ but it has no meaning and is
omitted.

The relation $\succ$ is a  partial order on $\tilde{V}.$ If
$\tilde{u}\succ\tilde{v}$ then all boundary sets $B_u$ contain
copies of all boundary sets $B_v,$ and $\dim B_u\ge\dim B_v$ for
all choices $u\in\tilde{u}, v\in\tilde{v}.$ Thus the largest
classes of boundary sets are the faces, then we have
one-dimensional sets, Cantor sets, and as terminal classes
countable and finite sets. Boundary sets from classes
$\tilde{u},\tilde{v}$ with $\tilde{u}\succ\tilde{v}$ can also have
equal dimensions, however.

\begin{Rem}
Consider a class $\tilde{u}\in \tilde{V},$ as a subset of $V.$
There are three cases.
\begin{enumerate} \item[ (i)] $\tilde{u}$
contains only one vertex. Then the boundary set $B_u$ has the same
dimension as the successor vertex $v$ (in case of several
successors, the union of the successor boundary sets). Namely
$B_u=f_i(B_v).$
\item[(ii)] $\tilde{u}$ is a cycle in $G,$ without loops and diagonals.
In that case the boundary sets $B_u$ with $u\in \tilde{u}$ are
countable unions of copies of the successor boundary sets.
\item[(iii)] $\tilde{u}$ contains two cycles. Only in this case $\tilde{u}$
produces its own Cantor structure, and the $B_u$ can have larger
dimension than the successor sets.
\end{enumerate}
\end{Rem}

In Figures \ref{ng6}, \ref{ng7}, and \ref{ng4} we have one class
of faces, a number of singleton classes and a terminal cycle which
is a class of two, one or four points. In Figure \ref{ng2} we have
a class of six faces, and a class of eight infinite neighbors h..r
on the right, as well as three singleton classes. The vertex s
describes a point neighbor, the vertices g and p correspond to
infinite boundary sets with equal dimensions as h..r, by the above
remark, (i). We now have to decide whether g,p and h..r are faces.

\section{Identifying faces: the general case}
{\bf An algorithmic approach to find the faces. } When we consider
one equivalence class of the neighbor graph which contains a face,
then it consists of faces only.

\begin{Prob} Must all faces belong to one equivalence class of $G$
?\end{Prob}

This problem also arises when we work with the contact matrix. In
\cite{al}, irreducibility was just assumed.  We present methods to
decide the problem algorithmically. We start with a simple
topological observation.

\begin{Prop} The faces cover the boundary of a self-affine lattice
tile $T.$
\end{Prop}

\begin{Proof} Let $x$ be a boundary point of $T$ and $U$ an open
neighborhood of $x.$ Let $E=U\cap\partial T$ and let $T_k,
k=1,...,q$ denote the neighbor tiles of $T.$ The sets $T\cap T_k$
cover $E,$ and by the Baire category theorem one of these sets,
say $T\cap T_1,$ must contain an open subset $V$ of $E.$ Each
$y\in V$ has a neighborhood $W$ in $\R^n$ which is contained in
$T\cup T_1.$ So $T\cap T_1$ is a face. Since this holds for any
$U,$ the point $x$ is an accumulation point of the faces. Since we
have only finitely many faces, $x$ must belong to one of them.
\end{Proof}

\begin{Thm} (The neighbor graph decides which boundary sets are faces)\\
Let us assume that all neighbors of the tile $T$ are compatible.
Then a boundary set $B_u$ is a face if it is not contained in the
union of all other boundary sets.\\
For the neighbor graph this means that there is a word $\ii$ which
belongs to the language $L_u$ of labels of paths starting in $u,$
but not to $L_v$ for any other $v\not= u.$
\end{Thm}

\begin{Proof} Consider the union of the tile $T$ and all its
neighbor tiles. If $B_u$ is not a face, then by the proposition it
is contained in the union of the remaining boundary sets.  If
$B_u$ is a face, there is $x\in B_u$ and a neighborhood $U\subset
T\cup T+u$ of $x.$ Then $U$ cannot intersect any other neighbor
$T+v$ since neighbors must have no common interior points. The
first part of the theorem is proved.\\
Now $x$ has an address: $x=\pi (\ss)$ where $\pi$ is the
projection from the symbol space $I^\infty$ to $T.$ Since $\pi$ is
continuous \cite{br,fa}, there is a prefix $\ii$ of $x$ such that
all addresses $\tt$ which begin with $\ii$ will fulfil $\pi
(\tt)\in U.$ By Theorem 4.5 (iii), this shows that for a face
$B_u$ there is a word $\ii$ which is contained in $L_u$ and in no
other $L_v.$\\
Conversely, suppose $B_u$ is not a face and $\ii\in L_u.$ There is
an address $\ss$ which begins with $\ii$ such that $\pi(\ss)=x$ is
not on the intersection of different pieces $f_j(T).$ (It is known
that the intersections have Lebesgue measure 0 \cite{bg}, and
since the normalized Lebesgue measure on $T$ can be considered as
the image measure of the product measure
$\{\frac{1}{m},...,\frac{1}{m}\}$ on $I^\infty,$ no cylinder can
be mapped completely into the intersections.) Since $x$ is also in
a face $B_v$ and has only one address, $L_v$ contains the sequence
$\ss$ starting with $\ii .$
\end{Proof}

\begin{figure}
 \begin{center}
  \includegraphics[width = 75mm]{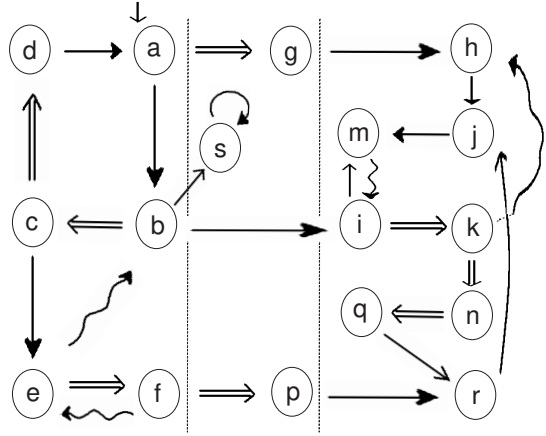}
  \caption{\label{ng2} Neighbor graph of
  twindragon \EE. The strong component of a..f represents
  faces, s is a point neighbor. Vertices h..r form a strong component
  representing infinite boundary sets. g and p are components with
  boundary sets of the same dimension as h..r.}
  \end{center}
\end{figure}

{\bf Examples of address calculations for boundary sets. } We
shall apply our theorem to Figure \ref{ng2}. Address calculations
are boring and best left to the computer, but they can provide
valuable information about the intersection of pieces. We start
with simple examples.

\begin{Exa} ( Point neighbors in \CC and \BB)\\
Which faces contain the point neighbors of \CC and $ \mathcal B ?$
In Figure \ref{ng6} we find
\[ B_{\rm g}=\pi(\overline{2})=B_\aa\cap B_\bb \cap B_\cc \cap
B_{\rm d}\, , \qquad  B_{\rm h}=\pi(\overline{2211})\subseteq
B_\cc\cap B_{\rm d}\cap B_{\rm e}\quad \mbox{ and} \]
\[ B_{\rm i}=\pi(\overline{2112})=B_{-\cc}\cap B_{\rm e}
\cap B_{\rm f}\, , \qquad  B_{\rm
j}=\pi(2\overline{1221})\subseteq B_\aa\cap B_{-\cc}\] so g,h,i
seem proper vertices but j seems to be on an edge. On the other
hand, both point neighbors of Figure \ref{ng7} coincide with the
intersection of only two faces:
\[ B_{\rm i}=\pi(\overline{2})=B_{\cc}\cap B_{\rm e}, \qquad
B_{\rm h}=\pi(1\overline{2})=B_{\aa}\cap B_{\cc}\ .
\]
\end{Exa}

Already these calculations show the strange polyhedral structure
of twindragons. For Figure \ref{ng2} the boundary sets look much
more complicated than for ordinary polyhedra where an edge is
covered by two faces. The proof of the next proposition will also
reveal hidden symmetry behind the apparently orderless structure
of $G.$

\begin{Prop} ( Boundary structure of \EE)\\
The vertices {\rm g, p, h..r} in the neighbor graph of the
twindragon \EE do not represent faces. The boundary set {\rm i} is
contained in the union of the four faces {\rm e, -f, c, -a} and
not in a union of three of them.
\end{Prop}

\begin{Proof}
We determine the language $L_{\rm i}$ from the three elementary
cycles which lead in $G$ from i to either i or -i (cf. section 7
and Example 9.4).
\begin{equation}\label{Li} L_{\rm i}= *2111 L_{\rm i} \cup
12 L_{\rm -i} \cup ***1122 L_{\rm -i}= C_1L_{\rm i}\cup C_2L_{\rm
i}\cup\bigcup_{j=0}^4 D_jL_{\rm -i}
\end{equation} where $*$ is considered as a wildcard for both 1
and 2, and $C_1=12111, C_2=22111, D_0=12, D_1=1*11122,
D_2=1*21122, D_3=2*11122, D_4=2*21122.$ We now consider the set of
faces $J=\{ {\rm e, -f, c, -a}\} .$ We will show
\begin{equation}\label{li}
L_{\rm i}\subset L_J=L_{\rm e}\cup L_{\rm -f}\cup L_{\rm c}\cup
L_{\rm -a}.
\end{equation}
First we study the action of $C_1$ and
$C_2$ on the vertices of $J.$ In the following table $+$ means
there is path from x to y with label $C_1,$ and $*$ means the same
for $C_2.$

\begin{center}\begin{tabular}{|l|c|c|c|c|}\hline
x\ $|$y&e&-f&c&-a\\ \hline e&$*$&&$*$&\\ \hline -f&&$+$&&$+$\\
\hline c&&$*$&&$*$\\ \hline -a&$+$&&$+$&\\ \hline \end{tabular}
\end{center}

{\it Claim 1: } any word $C\in \{ C_1,C_2\}^*$ is in $L_J,$ and it
can be realized with any prescribed terminal point in $J.$\\
This claim can be proved by induction on the number $t$ of terms
in $C=C_{q_1}^{n_1}...C_{q_t}^{n_t},$ using the table. The
starting point always depends on the word. For instance, it must
be -f for $C=C_1^3.$

Now we study how the $D_j, j=0,..,4$ lead from $J$ to $-J.$ In $G$
it can be seen that $D_0$ leads from e to -e, -f to f, c to a, and
-a to -c, respectively. Moreover, $D_1$ leads from -a to all
vertices of $-J,$ and from no other vertex of $J$ to any vertex of
$-J.$ Similarly, $D_2$ leads exclusively from -f to all of $-J,$
and $D_3$ leads exclusively from e to all of $-J,$ and $D_4$ leads
exclusively from c to all of $-J.$ From this observation it
follows that all vertices of $J$ are needed to obtain $L_{\rm
i}\subset L_J.$ Moreover, it is clear that

{\it Claim 2: } any word $D_j, j=0,..,4$ can be realized by a path
from $J$ to $-J$ with any prescribed terminal point in $-J,$
provided we can select the starting point.

Applying Claims 1 and 2 alternatingly, to $J$ and its opposite set
$-J,$ we complete the proof of (\ref{li}). The rest follows from
Theorem 9.3.
\end{Proof}

\section{Dimensions and eigenvalues}
{\bf The modified Hausdorff dimension. } In this section we give
an alternative proof for the fact that \EE has only 12 faces. It
is based on recent work of He and Lau \cite{hl} and of Akiyama and
Loridant \cite{al}. Let $H$ denote the adjacency matrix of the
graph $G\setminus \{ 0\}.$ For two non-zero vertices $u,v$ the
entry $h_{uv}$ denotes the number of edges leading from $u$ to
$v.$ It is also possible to take the adjacency matrix of our
simplified graphs, where only one vertex from each pair $\{
k,-k\}$ is taken into account and $h_{uv}$ counts the edges from
$u$ to $\{ v,-v\}.$ Since $H$ is a non-negative matrix it has a
real eigenvalue of maximum modulus which we call $\lambda .$

When the matrix $M$ defining the affine maps $f_j$ is conjugate to
a similarity matrix with expanding factor $R$ (see Proposition
2.2) it was shown by Mauldin and Williams \cite{mw} that the
Hausdorff dimension of the sets $B_k$ from (\ref{gradir}) fulfils
\begin{equation}\label{dimb}
\dim B_k \le \frac{\log\lambda}{\log R}\ ,
\end{equation}
and equality holds if an open set condition is satisfied, see the
discussion at the end of this section.

For self-affine sets, however, calculation of dimensions is much
more complicated then for self-similar ones. However, for our case
where all $f_j$ are defined by one matrix $M,$ He and Lau
\cite{hl} have defined a pseudo-norm $w$ on $\R^n$ which turns the
$f_j$ into similarity mappings, with factor $r=1/\sqrt[n]{\det
M}.$ Moreover, $w$ fulfils a weak type of triangle inequality
which is sufficient to develop Hausdorff measure and dimension as
usual. And $w$ generates the Euclidean topology even though it
drastically distorts the geometry of $\R^n.$ The value of the
modified Hausdorff dimension $\dim^w B_u$ was used in \cite{hl} to
give estimates of the real Hausdorff dimension in terms of
eigenvalues of $M.$ Akiyama and Loridant \cite{al} applied the
modified Hausdorff dimension to boundary sets of tiles. Here we
shall compare it with the topological dimension.

\begin{Thm} (Upper estimate of topological dimension of boundary sets)\\
Suppose that $T$ is a self-affine tile with $m$ pieces in $\R^n,$
and all neighbors are compatible. Let $W$ denote a strong
component of the graph $G,$ and let $\lambda_W$ denote the largest
eigenvalue of the adjacency matrix $H_W$ of $W$ as a subgraph of
$G.$ If for some integer $q<n,$
\[ \frac{n\log\lambda_W}{\log m}< q \]
then the topological dimension of all boundary sets $B_v$ with
$v\in W$ is strictly smaller than $q.$
\end{Thm}

\begin{Proof}
The matrix $H_W$ is irreducible, and $W$ is non-periodic (i.e. the
g.c.d. of the lengths of cycles from a point $v\in W$ to itself is
1) because of the compatibility of neighbors. So $\lambda_W$ is
the unique largest eigenvalue. The Hausdorff dimension of the sets
$B_v$ with $v\in W$ with respect to the He-Lau pseudo-norm is $\le
\frac{\log\lambda_W}{\log R}=\frac{n\log\lambda_W}{\log m}$ since
$R=\sqrt[n]{m}$ \cite{al,hl}. The pseudo-norm generates the
Euclidean topology, and the topological dimension is always
smaller than the Hausdorff dimension.
\end{Proof}

\begin{Exa} ( Boundary sets of the twindragon \EE)\\
For \EE we determine the dimension of the boundary sets $B_v$ with
$v\in\{\rm g, p, h..r\} .$ We obtain $\lambda_W \approx 1.554.$ A
small polynomial for $\lambda_W$ can be obtained from equation
(\ref{Li}) for $B=B_{\rm i}$ which expresses $B$ as a self-affine
set
\[ B= \bigcup_{i=1}^2  f_{i2111}(B)\cup
f_{12}(-B)\cup\bigcup_{i,j,k=1}^2 f_{ijk1122}(-B)\ . \] The
characteristic equation is $1=2z^5 +z^2+8z^7,$ which for
$\lambda=1/z$ becomes \  $\lambda^7-\lambda^5-2\lambda^2-8=0.$

The modified Hausdorff dimension is $\log\lambda_W
/\log\sqrt[3]{2} \approx 1.908.$ Thus the topological dimension of
the $B_v$ is $\le 1.$ So the $B_v$ cannot be faces. A calculation
for the faces {\rm a..f} gives modified dimension 2.67.
\end{Exa}

For the complicated twindragons \FF and \GG we also computed the
strong components and their modified dimensions. For \GG, the
above theorem does not apply since beside the irreducible
component of 22 faces with modified dimension 2.87, there is
another component of 20 boundary sets, plus 2 singleton
components, all with modified dimension 2.13. All the remaining 32
boundary sets have modified dimension 1.52. This boundary seems
really complicated! The other example is less intricate. We omit
the details of the following example.

\begin{Exa} ( Boundary sets of the twindragon \FF )\\
For \FF, the neighbor graph $G\setminus\{ 0\}$ has an irreducible
component $V_1$ of 16 faces with modified dimension 2.78, and the
other 32 neighbors are not faces. 10 of them form a cycle without
diagonals as second component $V_2,$ and there are two singleton
components between $V_1$ and $V_2.$ Since the cycle $V_2$ consists
of six double arrows and 4 simple arrows, $\lambda_2=2^{3/5}$ and
$3*\log\lambda_2/\log 2= 1.8.$ Due to the theorem, the topological
dimension is at most 1.

Moreover, there is another 10-cycle $V_3$ with only two double
arrows and $\lambda_3=2^{1/5},$ and 10 singleton components
between $V_2$ and $V_3.$ So the 20 corresponding boundary sets
have modified dimension $3*\log\lambda_3/\log 2= 0.6$ and by
Theorem 10.1 must be Cantor sets.

When one considers the labels along the cycles $V_2$ and $V_3,$
one obtains the languages
\[ L_2=\overline{*1**1*2**2}\quad\mbox{ and }\quad
L_3=\overline{11*1122*22} \] and their cyclic permutations as
address sets of the corresponding boundary sets. From this it is
easy to conclude that each Cantor set associated to a vertex in
$V_3$ is the intersection of two respective boundary sets in
$V_2.$\end{Exa}

We did not use the open set condition of the boundary equations
(\ref{gradir}), and this condition need not always be fulfilled.
In the cube (Example 6.3), for instance, some of the equations
(\ref{gradir}) contain equal terms since an edge is counted two
times. However, the open set condition will be true if we remove
such double sets and if the neighbors are all compatible. In
particular it is true when we restrict ourselves to faces
\cite{al}. When the graph is a cycle, as $V_2$ and $V_3$ in the
last example, the open set condition also holds.

\section{Polyhedral structure of tiles}
{\bf Homeomorphism with a polyhedron. } We have seen that the
structure of boundary sets can be quite complicated. We are
interested in the presence of simple structure, however. Let us
say that $T\subset\R^3$ {\it has the structure of polyhedron} if
$T$ is homeomorphic to a ball in $\R^3,$ the faces are
homeomorphic to disks in $\R^2,$ and their non-empty intersections
are either singletons or homeomorphic to an interval - in such a
way that $T$ together with its boundary structure becomes
homeomorphic to a polyhedron in Euclidean $\R^3.$\smallskip

This condition is stronger than just requiring that $T$ is
homeomorphic to a ball. However, this is the structure which
crystallographers expect of a tiling: ``We define a tiling as a
periodic subdivision of space into bounded, connected regions
without holes, which we call tiles. If two tiles meet along a
surface, we call the surface a face. If three or more faces meet
along a curve, we call the curve an edge. Finally, if at least
three edges meet at a point, we call that point a vertex''
\cite{fdhkm}.\smallskip

{\bf Definition of edges and vertices. } We follow this quotation
and define a polyhedral structure on every self-affine tile. We
remove from $G$ all vertices which do not represent faces. Then we
define the graph $G^2$ of edges - that is, of intersection sets
$B_k\cap B_{k'}.$ For these, it is not so difficult to check
whether they are homeomorphic to intervals \cite{bm}. If
necessary, we can also determine the graph $G^3$ of vertices of
the tiling. Euler's polyhedra formula can then be used to check
whether $T$ is homeomorphic to an ordinary polyhedron. The method
is quite fast in showing that certain tiles are {\it not}
ball-like.

A definition similar to the following was given by Thuswaldner and
Scheicher \cite{st} but it was not applied to examples.

\begin{Def}\label{defng2}
Let $T\subset \R^n$ be a self-affine lattice tile with neighbor
graph $G=(V,E=\bigcup_{i=1}^m E_i).$ For $\ell =2,3$ the graph of
$\ell$-intersections of boundary sets of $T$ is
\[ G^\ell =(V^\ell ,E^\ell ) \quad\mbox{ with }\quad
V^\ell =\left\{ K=\{k_1,...,k_\ell\}\, |\mbox{ all }k_i\mbox{
different, } B_{k_1}\cap ...\cap B_{k_\ell}\not=\emptyset\right\}
\]
\begin{equation}\mbox{ and }\qquad  E^\ell =\bigcup_{i=1}^m E^\ell_i\quad
\mbox{ with } \label{neigra2}\end{equation}
\[ E^\ell_i=\{ (K,K',i)\, |\mbox{
there is a 1-1-map }\phi:K\to K' \mbox{ with }(k,\phi(k),i)\in
E_i\mbox{ for }k\in K\} \]\end{Def}

The idea is that a neighbor intersection $(T+k)\cap (T+k')$ has a
piece which is at the boundary of $T_i$ if and only if both
neighbor sets have such a piece, and the pieces do intersect. The
algorithmic check for non-empty intersection is the same as for
$G.$ We start with all $\ell$-subsets of $V,$ and step by step we
remove those which have no outgoing vertices to the remaining
vertex set. Since for $\ell=2,3$ we must have non-empty
intersections, we end with a graph where each set has outgoing
vertices, and hence `infinite paths' of edges (counting repeated
use of cycles). The properties of $G^\ell$ are similar to those of
$G.$

\begin{Prop} $G^\ell$ provides a graph-directed fractal structure
on the $\ell$-intersec\-tions of neighbors. For each vertex $K\in
V^\ell,$ the labels of infinite paths in $G^\ell$ with starting
point $K$ are the addresses of the corresponding boundary
intersection $\bigcap_{k\in K} B_k .$ The cardinality of this
intersection can be determined as in Theorem 4.5.
\end{Prop}

\begin{Exa}
We continue Example 4.4 by determining $G^2$ for the six faces of
the Sierpi\'nski tile, see Figure \ref{siertile}. Since the
outgoing edges from $-\aa, -\bb, -\cc$ all start with $4,2,$ and
3, respectively, the corresponding boundary sets have pairwise
empty intersection. The intersection of a boundary set with its
opposite is also empty. We get edges
\[\{\aa, -\bb\}\stackrel{2}{\to}\{\aa, -\cc\}\stackrel{3}{\to}
\{\aa, -\bb\}\] so $B_\aa \cap B_{-\bb}$ is a singleton with
address $2323...=\overline{23}$ and $B_\aa \cap B_{-\cc}$ has
address $\overline{32}.$ The vertices $\aa,\bb,\cc$ all have a
loop with label 1, so the corresponding three boundary sets meet
in the point with address $\overline{1}.$ Altogether, $G^2$
consists of three cycles of length 2, and three isolated loops.
$G^3$ has only one vertex with a loop, representing the center
point of the tile.
\end{Exa}

\section{Polyhedral structure of twindragons}
{\bf Intersections between antipodal surfaces. } The first
information which the face intersection graph $G^2$ provides is a
list of non-empty intersections of pairs of faces. As in Theorem
4.5, it is easy to decide whether such an intersection consists
just of one point.

\begin{Prop}
Two faces $B_u, B_v$ with $u,v\in V$ intersect each other if and
only if $\{ u,v\}$ is a vertex in $G^2.$ Moreover, if there is
only one infinite directed path in $G^2$ starting in the vertex
$\{ u,v\} ,$ then $B_u\cap B_v$ is one point, and the address of
this point can be read from the labels of the path.
\end{Prop}

For our twindragons we check whether there is a face $B_k$ which
intersects its opposite face $B_{-k}.$ If such intersections
exist, it is unlikely that our tile has the structure of a
polyhedron. Actually, this happens for four of our twindragons and
provides an argument to show that \DD and \EE are not homeomorphic
to a ball.

A set $A\subset \R^n$ is called simply connected if it is
connected and its homotopy group is trivial. The latter means that
$A$ has ``no holes'' -- any closed curve $C\subset A$ can be
contracted within $A$ to a point.

\begin{Thm} (Tiles for which the interior is not simply
connected)\\
Let $T$ be a connected self-affine lattice tile which has a
symmetry center $c$ which is not an interior point of $T.$ Then
the interior of $T$ is not simply connected.\\
In particular, the interior of \DD and \EE is not simply
connected. These twindragons are therefore not homeomorphic to a
ball.
\end{Thm}

\begin{Proof} We can assume $c=0$ and thus $T=-T$ by changing the maps as in
Proposition 6.1. If 0 is not an interior point of $T,$ there is
some tile $T+k$ in the tiling which contains 0, and by symmetry
also $-(T+k)=T-k$ does not intersect $\inte T$ and contains 0.
Take a big sphere $S$ around 0 which contains the three tiles in
its interior. Let $z\in S$ be of the form $z=t\cdot k, t>0.$ There
is a line segment $\alpha$ from $z$ to the point $y\in T+k$ which
maximizes the scalar product $<y,k>.$ Since $T$ is connected, it
is arcwise connected, and there is an arc $\beta\subset T+k$ from
$y$ to 0. The union $\varepsilon$ of the arcs $\alpha , \beta ,
-\beta, -\alpha$ is an arc in $S$ from $z$ to $-z$ which does not
intersect ${\rm int}\, T.$

Now assume ${\rm int}\, T$ is connected and take a point $x\in
{\rm int}\, T.$ The opposite point $-x$ is also in ${\rm int}\,
T.$ So there is a connecting arc from $x$ to $-x.$ Covering this
arc by finitely many open balls inside ${\rm int}\, T,$ we see
that we can replace it by a union $\gamma$ of finitely many line
segments, and we can choose them so that no two segments are
parallel. Then the union $\gamma\cup -\gamma$ is a closed curve
from $x$ through $-x$ to $x$ with a finite number of
self-intersections. If $x_1$ is the first self-intersection point,
we find a closed arc from $x_1$ through $-x_1$ to $x_1$ with two
self-intersection points less. After a finite number of steps we
get a simple closed curve $\delta$ in ${\rm int}\, T$ from some
$x_n$ through $-x_n$ to $x_n.$ By construction, $\delta$ surrounds
the arc $\varepsilon$ and thus cannot be contracted within ${\rm
int}\, T$ to a single point. So ${\rm int}\, T$ is not simply
connected. The second assertion will be proved in the example
below.
\end{Proof}

\begin{Prop}
If an infinite path in the neighbor graph $G$ of a twindragon
starts in the root and contains no double arrow, it defines an
address of the center 0 of $T.$
\end{Prop}

\begin{Proof}
For a path of single arrows which addresses a point $x,$ the
opposite path is obtained by just interchanging labels 1 and 2. In
the proof of Proposition 6.1, it was shown that the address of the
opposite path corresponds to the point $-x.$  For paths starting
in the root of $G,$ Theorem 4.5 (i) says that both addresses
belong to the same point. $x=-x$ implies $x=0.$
\end{Proof}

\begin{Exa} ( Faces intersecting their opposite face in \DD and
\EE )\\
For \DD, Figure \ref{ng4} shows only one path from the root
without double arrows, given by the cycle {\rm acfd} which
describes an 8-cycle in the complete neighbor graph. The
associated address $1\overline{22211112}$ corresponds to the
center $0.$ As can be directly seen in Figure \ref{ng4}, this
address is also obtained from a path starting in vertex \cc, going
through the cycle {\rm cefg} and the opposite vertices. Thus the
center 0 belongs to the boundary set $B_\cc,$ hence not to the
interior of $T.$

For \EE, two root addresses of the center $0$ are visible in
Figure \ref{ng2}: $121\overline{2}$ from {\rm abs}, and
$122\overline{1221}$ from {\rm abim.} The second address is also
obtained when we start at {\rm -e} and go {\rm -ebcef-e-fef... }
Thus $0\in B_{\rm -e},$ the center is not in $\inte T.$
\end{Exa}

The addresses were found from a calculation of $G^2$ although $G$
is sufficient to check them. We determined $G^2$ also for the more
complicated twindragons and found that for \FF there are three
faces $B_k$ which intersect their opposite faces $B_{-k},$ not in
a point, but in a Cantor set with $\dim^w B_k\cap B_{-k}=0.6.$ It
is not clear whether any of the $B_k$ contains $c,$ so the above
argument fails, but it is obvious that $T$ has not the structure
of a polyhedron. In an ordinary polyhedron, given topologically as
a planar map on the sphere, at most one pair of opposite faces can
intersect. For \GG the situation is similar: we have even five
pairs of opposite faces $B_k$ which intersect each other, not in a
point, but in an uncountable set.\smallskip

{\bf The polyhedral structure of \BB and \CC. } Only the two
twindragons shown in Figures \ref{twin6} and \ref{twin7} can still
have the structure of a polyhedron. We shall prove that this is
not quite the case but we expect that

\begin{Conj} The twindragons \BB and \CC are homeomorphic to a
ball.\end{Conj}

It will not be possible to settle this question here, but we shall
establish the polyhedral structure which leads us to the
conjecture. We start with the list of non-empty intersections of
faces which was established by computing the graph $G^2.$
\smallskip

\begin{Exa} ( Polyhedral structure of \CC)\\
In \CC, the following faces have uncountable intersection:\\
{\rm a \emph{with} b,f,-c,-e \qquad b \emph{with} a,c,d,-e\qquad c
\emph{with} b,d,-a,-f\\  d \emph{with} b,c,e,f\qquad\quad  e
\emph{with} d,f,-a,-b \qquad  f \emph{with} a,d,e,-c\\}
Moreover, there are one-point intersections\\
{\rm a,-a \emph{with} d,-d, \qquad b,-b \emph{with} f,-f\quad
\emph{ and }\quad c,-c \emph{with} e,-e\\ } and two-point
intersections \aa  with \cc  and {\rm -a} with {\rm -c}. The two
points coincide, however, since the corresponding addresses are
identified by Theorem 4.5 (i).
\end{Exa}

\begin{Proof} The proof is a simple but lengthy calculation.
We sketch some facts which can be seen directly from inspection of
$G$ in Figure \ref{ng6}. The address $\overline{2}$ belongs to
$L_\aa, L_\bb, L_{\rm d}, L_{\rm f}$ and it is clear that it is
the only address in $L_a\cap L_{\rm d}$ as well as in $L_\bb\cap
L_{\rm f}.$ For the opposite faces we have address $\overline{1}.$
Moreover, $L_\cc\cap L_{\rm e}=\overline{2211}=L_{\rm d}\cap
L_{\rm -a}$ and $L_\cc\cap L_{\rm -e}=\overline{1221}=L_{\rm
f}\cap L_{\rm -b}.$ We can consider the languages restricted to
the vertices of faces, due to Proposition 9.2. Nevertheless, we
see that these addresses also appear for the point neighbors g, h
and -j, respectively, so they should really represent corner
points in the tiling by $T.$

It is obvious from $G$ that $L_{\rm f}\cap L_{\rm -f}=L_{\rm
f}\cap L_{\rm -e}=\emptyset .$ From this we conclude
\[ L_\bb\cap L_{\rm e}=2(L_{\rm d}\cap L_{\rm -c})=22(L_{\rm
f}\cap L_{\rm -f})\cup 22(L_{\rm f}\cap L_{\rm -e})=\emptyset .\]
Now we can determine
\[ L_\aa\cap L_\cc= *(L_\bb\cap L_{\rm e})\cup 2(L_\cc\cap L_{\rm
e})\cup 1(L_\bb\cap L_{\rm f})=2\overline{2211} \cup
1\overline{2}\, ,\] and these two addresses are equivalent since
they label opposite paths from the root of $G.$
\end{Proof}

\begin{figure}
 \begin{center}
  \includegraphics[width = 77mm]{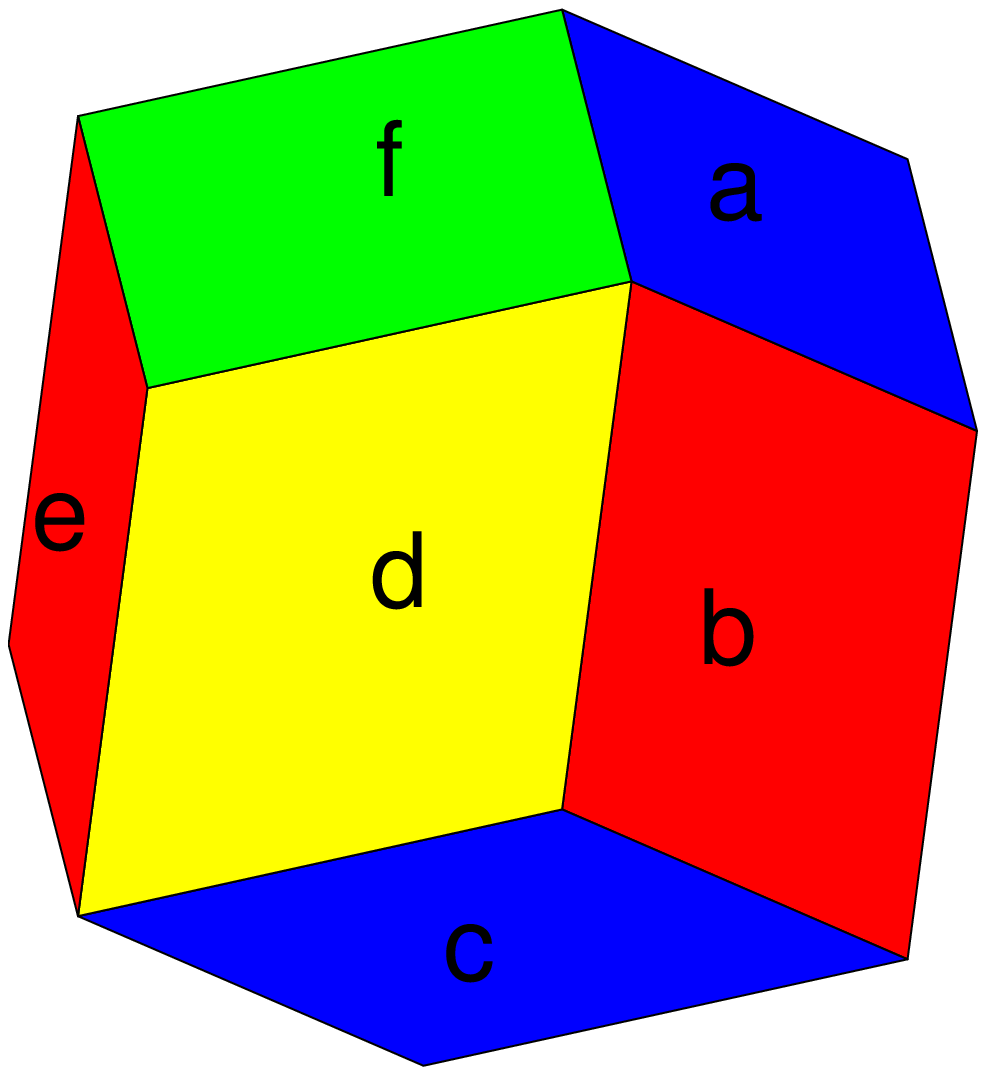}\quad
  \includegraphics[width = 67mm]{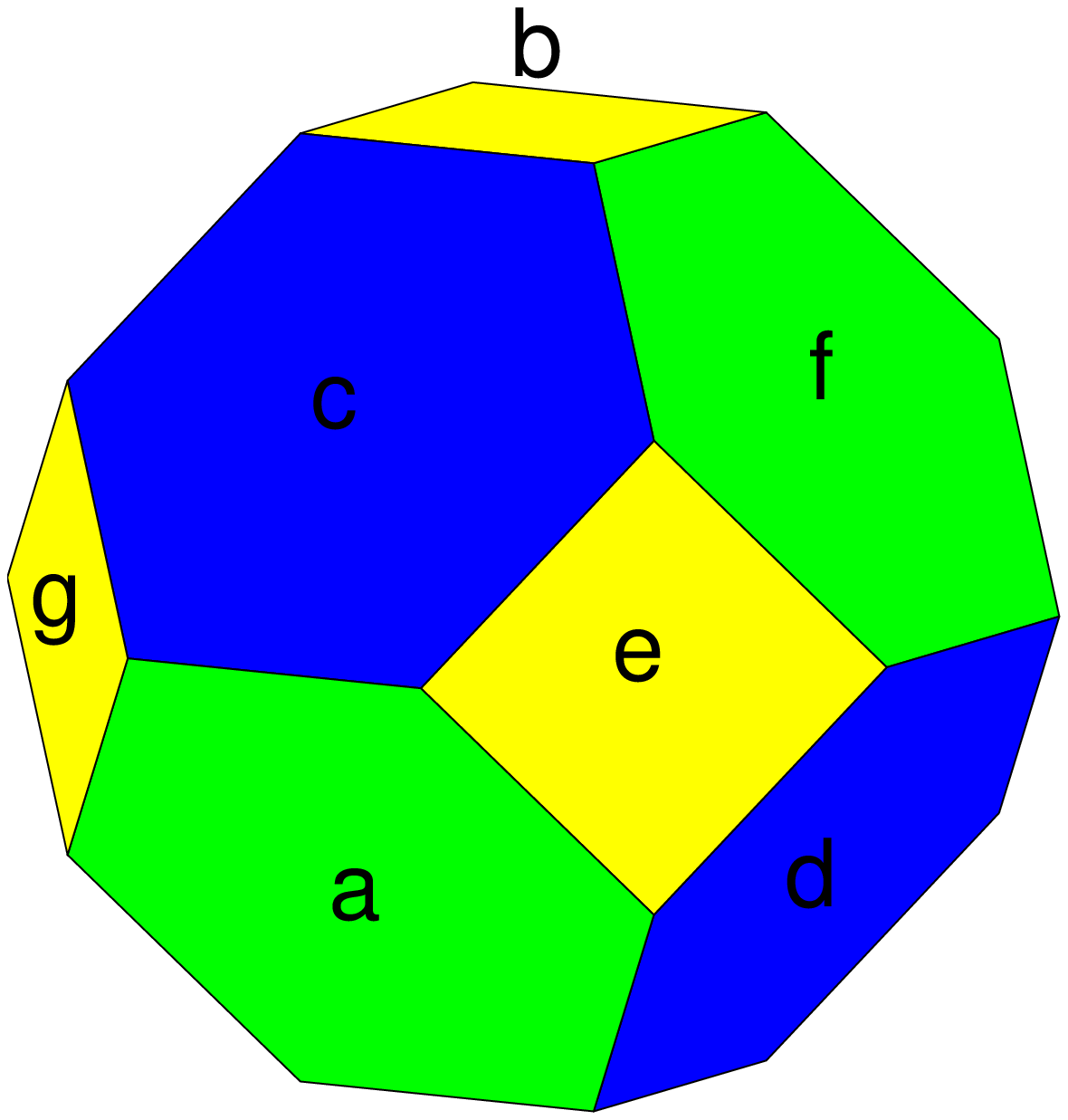}
  \caption{\label{twinpo} Polyhedral structure of
  \CC (left) and \BB (right)}
  \end{center}
\end{figure}

\begin{Prop} ( Equivalent convex polyhedron for \CC)\\
The polyhedron on the left of Figure \ref{twinpo} realizes exactly
all intersections of the faces of \CC, except for $B_\aa\cap
B_\cc$ and $B_{\rm -a}\cap B_{\rm -c}.$  It is an octahedron
truncated at four vertices, in the middle of the corresponding
edges, and its faces are rhombi.\\
Nevertheless, \CC has not the structure of an ordinary polyhedron
because the interior of the face $B_\bb$ is not connected and so
$B_\bb$ is not homeomorphic to a disk.
\end{Prop}

\begin{Proof}
The first part is done by checking all pairs of faces. For the
second part becomes clear when we prove that $B_\aa\cap B_\cc$
lies on the surface of $T,$ in $B_\bb .$ However, it can be seen
in $G$ that $1\overline{2}$ is an address in $L_\bb .$\\
Here is another argument. From $G$ we have the equation
$B_\bb=f_1(B_{\rm d})\cup f_2(B_{\rm d}).$ On the other hand
$f_1(T)\cap f_2(T)=f_1(B_\aa)=f_2(B_{\rm -a}).$ We have seen that
$B_\aa$ and $B_{\rm -a}$ intersect $B_{\rm d}$ in a single point,
and this also holds for their images under $f_1$ or $f_2.$ Thus
$f_1(B_{\rm d})$ and $f_2(B_{\rm d})$ have only one common point.
So $B_\bb$ consists of two isometric pieces which intersect in a
singleton, and cannot be homeomorphic to a disk.
\end{Proof}

It should be mentioned that the polyhedron in Figure \ref{twinpo}
describes \CC only topologically - the metric structure is
completely different. We have $B_{\rm f}=f_2(B_\aa)$ and $B_{\rm
d}=f_2(B_{\rm f})$ so these faces are of different size. The
intersection of the faces $B_\aa,B_\bb, B_{\rm d}, B_{\rm f}$ is
not a corner as in Figure \ref{twinpo}: the point $\pi
(\overline{2})=(\frac12,0,\frac12)'=\frac{\aa+{\rm
d}}{2}=\frac{\bb+{\rm f}}{2}$ is on the middle between the centers
of opposite adjoining faces (see the table at Figure \ref{ng6}). A
similar calculation shows that $B_\aa\cap B_\cc$ is in fact the
center point of $B_\bb .$  Now let us study \BB in the same way as
\CC .\smallskip

\begin{Exa} ( Polyhedral structure of \BB)\\
The following faces have uncountable intersection in \BB :\\
{\rm a \emph{with} c,d,e,g,-b,-f \qquad b \emph{with}
c,f,-a,-d\qquad c \emph{with} a,b,e,f,g,-d\\  d \emph{with}
a,e,f,-b,-c,-g \qquad e \emph{with} a,c,d,f \qquad  f \emph{with}
b,c,d,e,-a,-g\qquad g \emph{with} a,c,-d,-f.\\}
Moreover, there are one-point intersections\\
$B_\cc\cap B_{\rm -e}=\pi (\overline{12}),$ $B_\cc\cap B_{\rm
-a}=\pi (2\overline{21}),$ and their opposites.
\end{Exa}

\begin{Prop} ( Equivalent convex polyhedron for \BB)\\
The polyhedron on the right of Figure \ref{twinpo} realizes
exactly all intersections of the faces of \BB, except for the
point neighbors. It is an octahedron truncated at all vertices, at
one third of the corresponding
edges. The faces are regular hexagons and squares.\\
Nevertheless, \BB has not the structure of an ordinary polyhedron
because the interior of the faces $B_\bb$ and $B_{\rm d}$ is not
connected and so these faces are not homeomorphic to a disk.
\end{Prop}

This is proved similarly as for \CC . Note that faces of the
truncated octahedron never meet in a single point. Using Figure
\ref{ng7} we see that $B_\cc\cap B_{\rm -a}$ is on $B_\bb$ and
$B_\cc\cap B_{\rm -e}$ on $B_{\rm -d}.$

\section{Connectedness of the interior}
In this section, we go a small step towards proving our
conjecture, by proving that the interior of \CC is connected.

\begin{Prop}
Suppose $T=\bigcup_{j=1}^m f_j(T)$ is a self-affine tile and there
is a connected set $E\subset\inte T$ such that $E\cap
f_j(E)\not=\emptyset$ for $j\in I=\{ 1,...,m\}.$ Then the interior
of $T$ is connected.
\end{Prop}

\begin{Proof} $D_1=\bigcup_{j\in I} f_j(E)$ is a connected subset of
$\inte T.$ By induction we show that $D_n=\bigcup_{|w|\le n}
f_w(E)$ is connected, where $|w|$ denotes the length of the word
$w.$ Then $D=\bigcup_{n=1}^\infty D_n=\bigcup_{w\in I^*} f_w(E)$
is connected, and this set is dense in $\inte T.$ Hence $\inte T$
is connected.
\end{Proof}

\begin{Prop} $\inte$\CC is connected. \end{Prop}

\begin{Proof} A piece $T_w$ belongs to $\inte T$ if it is surrounded
in $T$ by all possible neighbors. $T_w$ has all possible neighbors
within $T$ if a path labelled with a suffix of $w$ leads from the
root to each vertex of $G$ (cf. proof of Theorem 5.1). When we
require this for faces only, $w=121212$ and $v=2121121$ fulfil
this condition, as one can check with Figure \ref{ng6}. We put
\[ E=T'_w\cup T'_v\cup T'_{\tilde{w}}\cup T'_{\tilde{v}}\, ,\qquad
\tilde{w}=212121, \tilde{v}=1212212\] where $T'_w$ means $T_w$
minus its one-point boundary sets. Then $E\subset \inte T.$ Since
$E$ contains $\pi (\overline{12})$ and $\pi (\overline{21}),$ the
set $E$ intersects $f_1(E)$ and $f_2(E).$ To apply Proposition
13.1, it only remains to verify that $E$ is connected. Since
twindragons are connected (cf. section 7), each $T'_u$ is
connected. Moreover, there are pairs of opposite paths from the
root in $G$ with labels $w1$ and $v,$ hence also with $\tilde{w}2$
and $\tilde{v},$ and also with $v$ and $\tilde{v},$ which all end
in a vertex of $G$ corresponding to a face. Thus $T_{w1}$ and
$T_v,$ $T_{\tilde{v}}$ and $T_{\tilde{w}2},$ as well as $T_v$ and
$T_{\tilde{v}}$ have a face in common. So $E$ is connected, and
the proposition applies.
\end{Proof}

$w$ and $v$ were obtained in a straightforward way. $w$ is the
shortest word for which $T_w$ has all possible face neighbors (we
used a computer search). $v$ is the shortest word with all face
neighbors which contains the center point 0 by Proposition 12.3
(this is easy: the address of 0 is $1\overline{212}$). In fact, 0
is an interior point of $T_v\cup T_{\tilde{v}}.$ For \BB  we also
determined such $v.$ The address of 0 is $1\overline{122},$ which
gives $v=1122122122.$ This indicates that the interior of \BB is
more fragmented than that of \CC and, to let the above method
work, we would need a longer chain of pieces.

\section{Remarks on algorithms}\label{algo}
{\bf Use of computer. } We have shown that the topology of
self-affine tiles can be studied by rather simple methods: finite
automata and regular languages. To answer concrete questions,
however, we have to perform quite a number of elementary logical
operations. Computer assistance seems necessary and convenient to
study more complicated examples.

Although all results in this paper have been checked by hand,
interactive computer work was essential to obtain them. Let us
briefly document the main algorithms we have used. \smallskip

{\bf Construction of the neighbor graph. } To construct $G,$ three
steps are performed, as indicated in Example 4.4.
\begin{enumerate}
\item[ (i)] Lower and upper bounds $l_q,b_q$ for $\{ x_q|\, x\in T\}$ are
derived for each coordinate $q=1,2,3.$ They need not be sharp, but
should be taken rather generously (we work with integers anyway).
If $b_q$ is too large, the computer will work a little longer --
if it is too small, only part of the neighbor graph will be
determined. For examples like \BB and \CC, 100000 points of an IFS
algorithm \cite{br} are completely sufficient to find $l_q$ and
$b_q$. For \FF and \GG there are rare outliers on the thin fibres,
and an exact estimate is needed.
\item[(ii)] Starting with $k=0,$ a list of new vectors is
calculated by the recursion \[ k'=Mk+k_j-k_i \] where $i,j$ take
all values in $\{ 1,..,m\}.$ The index of $k$ is listed together
with the labels $i,j$ at $k'$ as the predecessor of $k',$ and the
number of successors of $k$ is updated at $k.$ Of course, $k'$ is
only processed if it is within the bounds $l_q$ and $b_q.$
Moreover, it must always be checked whether $k'$ has already been
listed before. In that case, the new predecessor is added at the
old place. Since the number of possible $k$ is finite, the whole
list will be completed after finite time, and describes a graph
containing $G.$
\item[(iii)] All $k$ without successors are removed from
the list, and for each of their predecessors, the number of
successors is updated. This step is repeated until each $k$ has at
least one successor (or, in case of a Cantor set, $G$ is empty).
\end{enumerate}

When the resulting $G$ has predecessors of the root, the open set
condition is not fulfilled. This cannot happen if the $k_j$ form a
complete residue system for $M.$ \smallskip

{\bf Other operations with graphs. } For constructing the
intersection graph $G^2$ from $G,$ we have to consider $G\times
G,$ identify $(u,v)$ with $(v,u)$ and define the edges according
to Definition 11.1. Then we apply step (iii) as above, to exclude
empty intersections.

To identify point neighbors in $G,$ we use the adjacency matrix
$H$ defined in section 10. If $q$ is the number of non-zero
vertices in $G,$ then $k$ is a point neighbor if the sum of the
row associated with $k$ in $H,H^2,...,H^q$ always equals 1. When
we repeatedly apply a procedure like (iii) to the point neighbors,
we also remove the finite boundary sets. One-point intersection
sets in $G^2$ are found in the same way as point neighbors in $G.$

The partial order of boundary sets (section 8) is also determined
with $H.$ Let the matrix $N$ be the maximum of $H+H^2+...+H^q$ and
1, taken in each cell. Then $n_{uv}=1$ if and only if $u\succ v.$
In the row of $u$ there are the elements which can be reached from
$u,$ in the column of $u$ there are the $v$ from which $u$ can be
reached. The strong components of $G$ can be found by properly
ordering the rows and columns of $N.$

To deal with languages $L_u,$ we introduce adjacency matrices
$H_i$ for the edges with label $i.$ Then the matrix $H_1H_2H_2,$
for instance, tell us between which vertices there is a path
labelled $122.$ This fact was used to find the proof of
Proposition 9.5.

There are many similar tools waiting to be developed and tested to
help reveal the geometric structure of fractal tiles. \pagebreak

\begin{minipage}[t]{8.5cm}Christoph Bandt\\
Institute for Mathematics and Informatics\\
Arndt University\\ 17487 Greifswald, Germany\\
{bandt@uni-greifswald.de}\end{minipage}

\end{document}